\documentclass[12pt]{article}
\usepackage{setspace}

\usepackage{latexsym}
\usepackage{rotating}
\usepackage{amsmath}
\renewcommand{\epsilon}{\varepsilon}

\newcommand{\Var}{\mbox{Var}}

\usepackage{amsmath,bm}
\usepackage{amssymb}
\usepackage{ntheorem}
\usepackage[ansinew]{inputenc}
\usepackage{graphicx}
\usepackage{epsfig}
\usepackage{dsfont}
\usepackage{bbm}
\usepackage{subfig}
\usepackage{float}
\usepackage{url}
\usepackage{multirow}

\usepackage[authoryear]{natbib}
\newtheorem{satz}{Theorem}[section]
\newtheorem{example}{Example}[section]

\newtheorem{rem}[satz]{Remark}
\newtheorem{assumption}[satz]{Assumption}
\newtheorem{algo}[satz]{Algorithm}

\def\3{\ss}

\def \R{I \!\! R}   
\def \N{I \!\! N}   
\def \Z{Z \!\!\! Z} 
\def \C{\mathbb{C}}

\newcommand{\E}{\mathbbm{E}}

\newcommand{\bea}{\begin{eqnarray*}}
\newcommand{\eea}{\end{eqnarray*}}
\newcommand{\be}{\begin{eqnarray}}
\newcommand{\ee}{\end{eqnarray}}

\newcommand{\ba}{\begin{array}}
\newcommand{\ea}{\end{array}}
\newcommand{\cum}{\text{\rm cum}}

\newcommand{\Cov}{\text{\rm Cov}}
\def\3{\ss}

\begin{document}
\setstretch{1.635}

\title{Detection of multiple structural breaks in  multivariate time series}

\author{ Philip Preu\ss, Ruprecht Puchstein, Holger Dette \\
Ruhr-Universit\"at Bochum, 
Fakult\"at f\"ur Mathematik \\
44780 Bochum,   Germany 
}
\date{}
 \maketitle
\begin{abstract}
We propose a new nonparametric procedure for the detection and  estimation of  multiple structural breaks 
 in the autocovariance function of a multivariate (second-order) piecewise stationary process, which also identifies the
 components of the series where the breaks occur.  The new method is based on a comparison of the estimated spectral 
 distribution on different segments of the observed time series and  consists of three steps:  it starts with a consistent  
 test, which allows to prove the existence of structural breaks at a controlled type I error. Secondly,
it  estimates sets containing possible break points and finally these sets are reduced to identify  the relevant structural breaks 
and corresponding components which are responsible for the changes in the autocovariance structure. In contrast to all other methods  which have been proposed in the literature, our approach does not make any parametric assumptions, is
not especially designed for detecting one single change point and
addresses the problem  of multiple structural breaks  in the autocovariance function 
directly with no use of the binary segmentation algorithm.
 We prove that  the new procedure detects  all components and the corresponding
 locations where structural breaks occur with probability converging to one as the
  sample size increases and provide data-driven rules for  the selection of all regularization parameters. 
  The results are illustrated  by analyzing financial returns, and 
in  a simulation study it is demonstrated that the new procedure outperforms
 the currently available nonparametric methods for detecting breaks in the dependency structure of multivariate time series. 
\end{abstract}



\section{Introduction} \label{sec1}
\def\theequation{1.\arabic{equation}}
\setcounter{equation}{0}

The assumption of second order stationarity of time series is widely used in the statistical literature, because 
it   yields an elegant and powerful statistical methodology like parameter estimation or forecasting 
techniques. On the other hand  many real world phenomena can not be adequately described by stationary processes 
because these models are  usually not able to capture all stylized facts of an observed time series $(\boldsymbol{X_t})_{t \in \Z}$.
Consequently numerous alternative models have been proposed to address model features of non-stationarity [see \cite{jansen1981}, \cite{dahlhaus1997}, \cite{adak1998} or \cite{ombao2001} among many others], and one concept which became quite popular in the last decades is to assume that there is a certain number of segments on which the process can be considered as stationary [see for example \cite{starica2005} or \cite{fryzlewicz2006}]. This implies the existence of several points where the process 'jumps' from one stationary model to  another, and testing for the presence and   determining the location of such structural breaks is an important and challenging task. The main problems/questions  in this context are 
(1) Do there exist structural breaks?
(2)  If there exist  structural breaks, how many of them are present?
(3)  Where are the break points located?
(4)  In multivariate time series: in which components  do structural breaks  occur? \\ 
Because of its importance the problem of change point detection has found considerable interest in the statistical literature. 
For many decades, the identification of structural breaks  in the mean function was the dominating issue, and some early results are given by \cite{changemean}. 
Since this pioneering work numerous authors have considered this problem [see \cite{changemeanbiometrika}, \cite{banerjee1992}, \cite{bai1994}  or \cite{piotrjasa} among many others]. More recently  several authors have  argued that in applications the detection of 
changes in the autocovariance structure is of importance as well. Typical examples include the discrimination between stages of high and low asset volatility or the detection of changes in the parameters of an AR(p) model in order to obtain superior forecasting procedures. In the univariate case \cite{inclan1994} propose a nonparametric CUSUM-type test for changes in the variance of an independent
identically distributed sequence.  These results are generalized by \cite{leepark2001} to linear processes. 
\cite{changevariance} consider the same model and problem, but use the Schwarz information criterion. 
 \cite{aue2009} discuss multivariate time series and propose  a nonparametric test for structural breaks in the variance matrix, which is applicable to a broad class of stochastic processes.\\
On the other hand, the detection of changes in the complete autocovariance structure, that  is the existence of a change in the covariance function $\gamma(k)=\E(X_0X_k)$ at some lag $k \in \N$, is more complicated and  the corresponding literature is much less developed and mainly refers to parametric models. 
\cite{lavielle} use penalized methods  for identifying structural breaks and  \cite{lee2003} propose a CUSUM-type statistic for detecting changes in  specific parameters.  Recently \cite{davis2006} 
 suggest the minimum description length principle to fit piecewise  constant AR processes, while   \cite{motivation2} 
 use  AR-processes with structural breaks for modeling asset volatility.  \\ 
This list is of course far from being complete and could be extended by a vast amount of further articles on this subject. Nevertheless, to the best of our knowledge, there does not exist any nonparametric method for the detection of structural breaks in the autocovariance structure of multivariate time series which addresses all four problems raised at the beginning of this paper simultaneously. In fact, essentially all procedures either assume independence of the observed data or rely on parametric assumptions.
Moreover, these methods are particularly  designed to detect one structural break  and the application of the   binary segmentation algorithm 
[\cite{vostrikova}]  is used to locate multiple break points [see for example \cite{binaryFryzlevicz1} for a very recent approach]. 
Although this approach works both theoretically and practically, it is of course not necessarily optimal because it relies on procedures designed for detecting one change point. 

In this paper we present an integrative procedure for the detection and localization of structural breaks in the autocovariance structure 
of a multivariate time series, which addresses all raised issues simultaneously.  
Our approach  works
 in the frequency domain and is based on a comparison of the  spectral density on two consecutive blocks of  length, say $N$, 
  which is 'small' compared to the whole sample size, say $T$. 
 Roughly speaking, the method consists of three basic steps.   We start with the development of a new consistent bootstrap  
 test, which allows to prove the existence of structural breaks at a controlled type I error. 
 In this case, a set of possible break points is estimated, and finally this  set  is reduced to identify  the relevant structural breaks 
and corresponding components which are responsible for  these breaks.  The different steps are carefully described and illustrated in Sections \ref{sec3} and \ref{sec4}.
Here we also prove that the initial test keeps its nominal level, is consistent and that the detection rule identifies all components and the corresponding
 locations where structural breaks occur with probability converging to one as the
  sample size increases. \\
 As all methods in this context [see for example \cite{lee2003}]
the new  procedure depends on the choice of several regularization parameters. Therefore 
we develop  data-driven rules for  the selection of these parameters in Section \ref{sec5}, which can also be justified  theoretically and 
yield a good  performance of the procedure from a practical point of view. In Section \ref{sec6} we  investigate the  finite sample 
properties of the new method by means of a simulation study and  illustrate its application  analyzing
multivariate financial market returns. In particular  it
is demonstrated that in most cases the new proposal improves competing nonparametric methods for detecting multiple 
structural  breaks in the dependency structure of a multivariate time series.
Finally all  proofs are given in Section \ref{sec7}, in which some technical details are deferred to an additional Section \ref{sec8}. 

We conclude this introduction emphasizing that  the focus of this paper is on \textbf{nonparametric} detection of changes in the autocovariance structure of \textbf{multivariate} time series and \textbf{identification of
 the components} which are responsible
 for the structural breaks.  Nearly all articles on nonparametric change point detection in the dependency structure are restricted to the univariate case. 
 To the best of our knowledge,  only \cite{aue2009} provide a formal level-$\alpha$ test for the presence of structural breaks in a multivariate time series setting.  However, these authors do not investigate change point problems in the autocovariance structure, but consider the problem of testing for a constant covariance matrix. As a consequence their test is  not able to discriminate between changes in the dependency structure which keeps the variance constant. This situation occurs for example in an univariate AR(1) model where the AR parameter switches from $-a$ to $a$ for some $a \in (0,1)$. The procedure proposed in this paper solves problems of this type and in the multivariate case additionally  identifies the components, where the changes occur.

\section{Piecewise stationary processes} \label{sec2}
\def\theequation{2.\arabic{equation}}
\setcounter{equation}{0}

We assume to observe realizations of a centered $\mathbb{R}^d$ valued stochastic process $(\boldsymbol{X}_{t,T})_{t=1,...,T}$, where $\boldsymbol{X}_{t,T}=(X_{t,T,1},...,X_{t,T,d})^T$ has a piecewise stationary representation. This means that there exist an unknown number $K \in \N_0$ and points $0=b_0<b_1< \cdots <b_k<b_{K+1}=1$ such that
\vspace{-.75cm}
\begin{align}\label{defX_t}
\boldsymbol{X}_{t,T}=\sum_{l=0}^\infty \boldsymbol{\Psi}_{l}(t/T)  \boldsymbol{Z}_{t-l} \quad t=1,...,T,
\end{align}
where the functions $\boldsymbol{\Psi}_{l}:[0,1] \rightarrow \R^{d\times d}$, $l\in\mathbb{Z}$ are defined as $\boldsymbol{\Psi}_{l}(u)=\sum_{j=0}^K\boldsymbol{\Psi}_{l}^{(j)}1_{S_j}(u)$ and $1_{S_j}$ denotes the indicator function of the set $S_j=\{u:b_j  < u \leq   b_{j+1} \}$, $\{\boldsymbol{Z}_t\}_{t\in\mathbb{Z}}$ denotes a centered Gaussian White Noise process with covariance matrix $\boldsymbol{I}_d$ and the matrices $\boldsymbol{\Psi}_{l}^{(j)}\in \R^{d \times d}$ correspond to the piecewise constant coefficents of the linear representations on the segment $(\lfloor b_jT\rfloor,\lfloor b_{j+1}T\rfloor]$.
Throughout this paper we assume that $K$ is 'minimal' in the sense that for every pair $(i,i+1)$ with  $i\in \{0,...,K-1\}$  there exists an integer $l \in \N$ such that
$\boldsymbol{\Psi}_{l}^{(i)} \not= \boldsymbol{\Psi}_{l}^{(i+1)} $. This ensures that, if $K$ equals zero, there is no change point in the dependency structure, while structural breaks exist for $K \geq 1$. 
We also note that the assumption of Gaussianity is only imposed here to simplify technical arguments and the extension of the proposed methodology and results to more general innovations is straigthforward; see Remark \ref{gaussian} for more details.\\
We introduce 
$
\boldsymbol{f}_j(\lambda)=\frac{1}{2\pi}\sum_{l,m=0}^\infty\boldsymbol{\Psi}_{l}^{(j)}\big(\boldsymbol{\Psi}_{m}^{(j)}\big)^T \exp(-i\lambda(l-m))
$
and obtain for the $\C^{d \times d}$ valued time-varying (piecewise constant) spectral density matrix
\begin{equation}
\boldsymbol{f}(u,\lambda)=\frac{1}{2\pi}\sum_{l,m=0}^\infty\boldsymbol{\Psi}_{l}(u)\big(\boldsymbol{\Psi}_{m}(u)\big)^T \exp(-i\lambda(l-m))=\sum_{j=0}^K\boldsymbol{f}_j(\lambda)1_{S_j}(u).
\end{equation}
From this representation it follows that the spectral density has points of discontinuity in $u$ direction at the break points $b_i$  ($i=1,...,K$) whenever $K\geq1$. Therefore we propose to compare the spectral density $\lambda \mapsto \frac{1}{e}\int_{v-e}^v\boldsymbol{f}(u,\lambda)du$ with $\lambda \mapsto \frac{1}{e}\int_v^{v+e}\boldsymbol{f}(u,\lambda)du$ for some 'small' constant $e$. If there exist structural breaks,  the difference $\sup\limits_{\omega\in[0,1]}\frac{1}{e}|\int_0^{\omega\pi}\int_v^{v+e}\boldsymbol{f}(u,\lambda)dud\lambda-\int_0^{\omega\pi}\int_{v-e}^{v}\boldsymbol{f}(u,\lambda)dud\lambda|$ will be positive for $v \in \{b_1,...,b_K\}$ while it vanishes for  \linebreak{$v \in [0,1] \backslash\{b_1,...,b_K\}$} as $e\rightarrow 0$. In order to obtain a global measure for the presence of structural breaks, we define $||A||_{\infty}=\max\limits_{a,b=1,...,d}|A_{a,b}|$ as the usual maximum norm of a matrix $A\in\C^{d\times d}$ and consider
\vspace{-0.5cm}
\begin{equation}
D:=\sup_{v,\omega \in [0,1]}||\boldsymbol{D}(v,\omega)||_{\infty},\label{Dtheo}
\end{equation}
where for $v\in [e,1-e]$ and $\omega\in[0,1]$ the matrix $\boldsymbol{D}(v,\omega)$ is defined by
\begin{equation}
\boldsymbol{D}(v,\omega):=\frac{1}{e}\left(\int_{0}^{\omega\pi}\int_{v}^{v+e}\boldsymbol{f}(u,\lambda)dud\lambda-\int_{0}^{\omega\pi}\int_{v-e}^v\boldsymbol{f}(u,\lambda)dud\lambda\right) \in\mathbb{R}^{d\times d}\label{limit}
\end{equation}
 and we set $\boldsymbol{D}(v,\omega)=\boldsymbol{D}(e,\omega)$ and $\boldsymbol{D}(v,\omega)=\boldsymbol{D}(1-e,\omega)$ whenever $v\leq e$  amd $v\geq 1-e$ respectively.
Under the hypothesis of no structural break, i.e. $K=0$, we have $D=0$, while $D$ is strictly positive if structural breaks occur. In order to obtain a test for the null hypothesis
\vspace{-1cm}
\be \label{null}
H_0: \quad K=0,
\ee
it is therefore natural to estimate $D$ and to reject the null hypothesis for 'large' values of the estimator. We will construct such an empirical version of $D$ in the following section and derive its asymptotic properties. The identification of the location and components corresponding to the break points is based on an estimator of the components   $\sup_{\omega \in [0,1]}|[\boldsymbol{D}(v,\omega)]_{a,b}|$ and illustrated in Section \ref{sec4}.

\section{Testing for structural breaks} \label{sec3}
\def\theequation{3.\arabic{equation}}
\setcounter{equation}{0}

The first step of the proposed procedure consists in a statistical test which allows to prove the existence of structural breaks at a controlled type I error and is based on an empirical version of the quantity $D$ defined in \eqref{Dtheo}. To be precise we choose an even integer $N\leq T/2$ and consider the local periodogram 
\vspace{-0.5cm}
\begin{align}\label{locper}
 \boldsymbol{I}_N(u,\lambda):=\frac{1}{2 \pi N} \sum_{r,s=0}^{N-1} \boldsymbol{X}_{\lfloor uT \rfloor-   N/2   +1+s,T}\boldsymbol{X}_{\lfloor uT \rfloor-   N/2   +1+r,T}^T \exp(-i \lambda (s-r)) ,
\end{align}
where we use the convention $\boldsymbol{X}_{j,T} =0$ whenever $j \not \in \{1,\ldots , T\}$. Note that \eqref{locper} is the usual periodogram computed from the $N$ observations $\boldsymbol{X}_{\lfloor uT \rfloor-   N/2   +1,T}, \ldots, \boldsymbol{X}_{\lfloor uT \rfloor+   N/2   ,T}$, and it can be shown that the quantity $\boldsymbol I_N(u,\lambda)$ is an asymptotically unbiased estimator for the spectral density if $N \rightarrow \infty$ [see \cite{dahlhaus1997}].
An estimator of the matrix $\boldsymbol{D}(v,\omega)$ is then defined by
\vspace{-0.5cm}
\begin{equation}
\label{empprocess}
\hat {\boldsymbol D}_T(v,\omega):=\frac{1}{N}\sum_{k=1}^{\lfloor \omega   N/2   \rfloor}\Big(\boldsymbol I_N\big(v+N/(2T),\lambda_k \big)-\boldsymbol I_N\big(v-N/(2T),\lambda_k\big)\Big),
\end{equation}
if $v\in[\frac{N}{T},1-\frac{N}{T}]$ where $\lambda_k=2\pi k/N$ denote the Fourier frequencies. On the intervals $[0,\frac{N}{T})$ and $(1-\frac{N}{T},1]$ we define $\hat {\boldsymbol D}_T(v,\omega)$ as $\hat {\boldsymbol D}_T(\frac{N}{T},\omega)$ and $\hat {\boldsymbol D}_T(1-\frac{N}{T},\omega)$ respectively. So roughly speaking we construct an estimator of $\boldsymbol{D}(v,\omega)$  by replacing the integral by a Riemann sum, where the averaged time varying spectral density matrices $\frac{1}{e}\int_{v}^{v+e}\boldsymbol{f}(u,\lambda)du$ and $\frac{1}{e}\int_{v-e}^{v}\boldsymbol{f}(u,\lambda)du$ on the intervals $[v,v+e]$ and $[v-e,v]$ are replaced by the local periodograms $\boldsymbol{I}_N(v+N/(2T),\lambda)$ and $\boldsymbol{I}_N(v-N/(2T),\lambda)$. The final estimate of the quantity $D$ in \eqref{Dtheo} is then defined by 
\vspace{-.5cm} 
\begin{equation}
\label{supstat}
\hat D_T:=\sup\limits_{(v,\omega) \in[0,1]^2}||\hat {\boldsymbol D}_T(v,\omega)||_{\infty} =\max_{v \in [N/T,1-N/T]}\sup\limits_{\omega \in [0,1]}||\hat{\boldsymbol{D}}_T(v,\omega)||_{\infty}.
\end{equation}
The following results specify the asymptotic properties of the process \linebreak$\{\hat{\boldsymbol{ D}}_T(v,\omega)\}_{(v,\omega)\in[0,1]^2}$ under the null hypothesis \eqref{null} of no structural breaks and the alternative $H_1:K>0$ for different choices of the sequence $N$. Throughout this paper  the symbol $\Rightarrow$ denotes weak convergence in $L^\infty([0,1]^2)$ and we distinguish the cases $N/T\rightarrow 1/c\in(0,1/2)$ 
(Theorem \ref{hauptsatz}) and $N/T\rightarrow 0$ (Theorem \ref{hauptsatz2}). 
%
\begin{satz}\label{hauptsatz}
Suppose that the coefficients in the representation \eqref{defX_t} satisfy
\vspace{-.5cm} 
\begin{equation}
\sum_{l= 0}^{\infty}\sup_{u \in[0,1]}\| \boldsymbol{\Psi}_{l}(u)\|_\infty|l|<\infty\label{summ1},
\end{equation}
and that $N/T\rightarrow 1/c$ for some  $c\geq 2/\min_{i=1,...,K+1}|b_i-b_{i+1}|$ as $T\rightarrow\infty$. Then the following statements hold:
\begin{itemize}
\item[a)] If $K=0$,  the process $\{\sqrt{N}\hat{\boldsymbol{ D}}_T(v,\omega)\}_{(v,\omega)\in[0,1]^2}$ converges weakly to a centered Gaussian process $\{\boldsymbol{G}(v,\omega)\}_{(v,\omega)\in[0,1]^2}$, i.e.
\vspace{-0.5cm}
\begin{align}\label{austheorem}
\{\sqrt{N}\hat{\boldsymbol{ D}}_T(v,\omega)\}_{(v,\omega)\in[0,1]^2}\Rightarrow\{\boldsymbol{G}(v,\omega)\}_{(v,\omega)\in[0,1]^2}.
\end{align}
Here for $v_i\in[\frac{1}{c},1-\frac{1}{c}]$ $(i=1,2)$, the covariance kernel $Cov([\boldsymbol{G}(v_1,\omega_1)]_{a_1,b_1},[\boldsymbol{G}(v_2,\omega_2)]_{a_2,b_2})$
of $\{\boldsymbol{G}(v,\omega)\}_{(v,\omega)\in[0,1]^2}$ is given by
\vspace{-0.5cm}
\begin{align}
\begin{cases}
0\quad&\text{if }\quad\frac{2}{c}\leq |v_2-v_1|\\
-[2-|v_2-v_1|c]\frac{1}{\pi}\int_0^{\min(\omega_1,\omega_2)\pi}\rho_{a_1,a_2,b_1,b_2}(\lambda)d\lambda\quad&\text{if } \quad \frac{1}{c}\leq |v_2-v_1|\leq \frac{2}{c}\\
[2-3|v_2-v_1|c]\frac{1}{\pi}\int_0^{\min(\omega_1,\omega_2)\pi}\rho_{a_1,a_2,b_1,b_2}(\lambda)d\lambda \quad&\text{if }\quad 0\leq |v_2-v_1| \leq  \frac{1}{c}\\
\end{cases},\label{covariancekernel}
\end{align}
where $\rho_{a_1,a_2,b_1,b_2}(\lambda):=f_{a_1,a_2}(\lambda)f_{b_1,b_2}(-\lambda)+f_{a_1,b_2}(\lambda)f_{b_1,a_2}(-\lambda)$. If $v_i\notin[\frac{1}{c},1-\frac{1}{c}]$ for at least one $i\in\{1,2\}$ the covariance kernel is given by
\vspace{-0.5cm}
\begin{align*}
Cov([\boldsymbol{G}(v_1,\omega_1)]_{a_1,b_1},[\boldsymbol{G}(v_2,\omega_2)]_{a_2,b_2})=Cov([\boldsymbol{G}(a_c(v_1),\omega_1)]_{a_1,b_1},[\boldsymbol{G}(a_c(v_2),\omega_2)]_{a_2,b_2})
\end{align*}
where $a_c(v):=\min(\max(v,\frac{1}{c}),1-\frac{1}{c})$.
\item[b)] If $K \geq 1$,  there exists a constant $C \in \R^+$  with  $
\lim\limits_{T\rightarrow\infty}\mathbb{P}(\sup\limits_{\omega\in[0,1]}\big|[\hat{\boldsymbol D}_T(b_r,\omega)]_{a,b}\big|>C)=1 $
for all $(r,a,b) \in\{1,...,K\}\times\{1,...,d\}^2$ such that  $\sup\limits_{\omega\in[0,1]}|[{\boldsymbol D}(b_r,\omega)]_{a,b}|>0$.
\end{itemize}
\end{satz}

\begin{satz}\label{hauptsatz2}
Suppose that the coefficients in the representation \eqref{defX_t} satisfy \eqref{summ1} and that $N \rightarrow \infty$, $T^{\epsilon}/N\rightarrow 0$ and $N/T\rightarrow 0$ 
for some $\epsilon>0$ as $T \rightarrow \infty$.
Then the following statements are correct:
\begin{itemize}
\item[a)] If $K=0$, then we have for any $0<\gamma<\frac{1}{2}$: $\hat D_T=o_P(N^{-\gamma}).$
\item[b)] If $K \geq 1$, then part b) of Theorem \ref{hauptsatz}
holds.
\end{itemize}
\end{satz}

From these results it follows that under the null hypothesis of no structural breaks the expression $\hat D_T$ is of order $O_P(N^{-\gamma})$ where $\gamma=1/2$ (Theorem \ref{hauptsatz}) or $\gamma\in(0,1/2)$ (Theorem \ref{hauptsatz2}) while it is larger than some positive constant under the alternative. Note that in the situation of Theorem \ref{hauptsatz2}, i.e. $N/T\rightarrow 0$, it can also be shown that under the null hypothesis of no structural breaks the random variable $\sqrt{N}\hat{\boldsymbol D}_T(v,\omega)$ converges weakly to a Gaussian limit for any fixed pair $(v, \omega)$. However, in this case $\sqrt{N}\hat {\boldsymbol D}_T(v_1,\omega)$ and $\sqrt{N}\hat {\boldsymbol D}_T(v_2,\omega)$ are asymptotically uncorrelated whenever $v_1 \not= v_2$,  and therefore measurability is not given in the limit. 

Since the distribution under the null hypothesis depends on unknown quantities of the data generating process we propose resampling methods to obtain its quantiles. More precisely we will now develop a bootstrap procedure, which is closely related to the one dimensional AR($\infty)$ bootstrap introduced by \cite{Kreiss1988}. This methodology has found considerable attention in the recent literature [see \cite{choi2000}, \cite{goncalves2007} or \cite{bergpappolitis2010} among others] since it is easy to implement but has sufficient complexity to capture the predominant dependencies in the underlying process. In the present context it will yield critical values such that a test for structural breaks based on the statistic $\hat D_T$ is directly implementable. For a description and the statement of the theoretical properties of the resampling procedure we require the following central assumption.
\begin{assumption}\label{annahmenbootstrap}
The stationary $\mathbb{R}^d$-valued process $\{\boldsymbol{X}_t\}_{t\in\mathbb{Z}}$ with spectral density function $\boldsymbol{g}(\lambda)=\int_0^1 \boldsymbol{f}(u,\lambda) du$  has an AR($\infty$)-representation of the form
\vspace{-.5cm}
\begin{equation}\label{Ar1}
\boldsymbol{X}_t=\sum_{j=1}^{\infty}\boldsymbol{a}_j\boldsymbol{X}_{t-j}+\boldsymbol{\Sigma}^{1/2}\boldsymbol{ Z}_t,
\end{equation}
where $\{\boldsymbol{Z}_t\}_{t\in\mathbb{Z}}$ denotes a sequence of independent  $d$-dimensional $\mathcal{N}(0,\boldsymbol{I}_d)$ distributed random variables, $\boldsymbol{\Sigma} \in \R^{d \times d}$ is positive definite and $(\boldsymbol{a}_j)_{j\in\mathbb{N}}$ is a sequence of $d\times d$ matrices satisfying  
\vspace{-.85cm}
\bea
\det\Big(\boldsymbol{I}_d-\sum_{j=1}^{\infty}z^j\boldsymbol{a_j}\Big)\neq 0 \text{ for } |z|\leq 1\quad \text{and }\quad \sum_{j=0}^{\infty}|j|\|\boldsymbol{{a}}_{j}\|_\infty<\infty.
\eea
\end{assumption}
 The main motivation of the resampling procedure consists in the fact, that every stationary process can be approximated by an AR($p$) model if the order $p$ of the autoregressive process is sufficiently large. Therefore we choose an increasing sequence $p=p(T)\rightarrow\infty$ as $T\rightarrow\infty$ and approximate the process defined in \eqref{Ar1} by an AR($p$) model with coefficents
\begin{equation}\label{arpfitt}
(\boldsymbol{a}_{1,p},...,\boldsymbol{a}_{p,p}):=\underset{\boldsymbol{b}_{1,p},\ldots,\boldsymbol{b}_{p,p}}{\operatorname{argmin}} \text{tr} \Big(\mathbb{E}[(\boldsymbol{X}_{t}-\sum_{j=1}^p\boldsymbol{b}_{j,p}\boldsymbol{X}_{t-j})(\boldsymbol{X}_{t}-\sum_{j=1}^p\boldsymbol{b}_{j,p}\boldsymbol{X}_{t-j})^T]\Big)
\end{equation}
and innovations with  covariance matrix
$
\boldsymbol{\Sigma}_p=\mathbb{E}[ (\boldsymbol{X}_t-\sum_{j=1}^p\boldsymbol{a}_{j,p}\boldsymbol{X}_{t-j})(\boldsymbol{X}_t-\sum_{j=1}^p\boldsymbol{a}_{j,p}\boldsymbol{X}_{t-j})^T].
$
To be precise let $(\boldsymbol{\hat{a}}_{1,p},...,\boldsymbol{\hat{a}}_{p,p})$ denote a 'consistent' estimator of the minimizer in \eqref{arpfitt} (the precise assumptions regarding this estimator are specified in Theorems \ref{theorembootstrap} and \ref{theorembootstrap2} below and are, for example, fulfilled for the Yule Walker estimators). The bootstrap replicates $\hat D_T^*$ of $\hat D_T$ are then generated as follows: 
\begin{algo}\label{algo1} {\rm {\bf (autoregressive bootstrap)}  We  simulate data from the model
\vspace{-.5cm}
\begin{equation}
\label{fittedarp}
\boldsymbol{X}^*_{t,T}=\sum_{j=1}^{p}\boldsymbol{\hat{a}}_{j,p}\boldsymbol{X}^*_{t-j,T}+\boldsymbol{\hat{\Sigma}}_{p}^{1/2}\boldsymbol{Z}^*_{j},
\end{equation}
where $ \boldsymbol{\hat{\Sigma}}_{p}=\frac{1}{T-p}\sum_{j=p+1}^T(\boldsymbol{\hat{z}}_j-\boldsymbol{\bar{z}}_T)(\boldsymbol{\hat{z}}_j-\boldsymbol{\bar{z}}_T)^T,
$
the random variables $\boldsymbol{Z}^*_{j}$ are independent $\mathcal{N}(0,\mathbf{I}_d)$ distributed,
$\boldsymbol{\hat{z}}_{j}:=\boldsymbol{X}_{j,T}-\sum_{i=1}^{p}\boldsymbol{\hat{a}}_{i,p}\boldsymbol{X}_{j-i,T}$ $(j=p+1,...,T)$ and  $\boldsymbol{\bar{z}}_T:=\frac{1}{T-p}\sum_{j=p+1}^T\boldsymbol{\hat{z}}_j$.
The bootstrap statistic  $\hat D_T^*$  is defined  as the statistic $\hat D_T$ in \eqref{supstat}, where the observations $\{\boldsymbol{X}_{t,T}\}_{t=1,...,T}$ are replaced by its bootstrap replicates $\{\boldsymbol{X}_{t,T}^*\}_{t=1,...,T}$.
}
\end{algo}

To motivate this procedure note that under the null hypothesis of stationarity the process $\{\boldsymbol X_t\}_{t=1,...,T}$ and $\{\boldsymbol X_{t,T}\}_{t=1,...,T}$ have the same spectral density $\boldsymbol{g}(\lambda)=\int_0^1\boldsymbol{f}(u,\lambda)du=\boldsymbol{f}(\lambda)$. Now $\{\boldsymbol X_t\}_{t=1,...,T}$ is approximated by an AR($p$) process and $\{\boldsymbol{X}_{t,T}^*\}_{t=1,...,T}$ mimics this approximation, which follows from definition \eqref{arpfitt} and the consistency of the estimators $\hat{\boldsymbol{a}}_{j,p}$. Therefore the process $\{\boldsymbol{X}_{t,T}^*\}_{t=1,...,T}$ exhibits similar spectral properties as the stationary process $\{\boldsymbol X_t\}_{t=1,...,T}$. Consequently, under the null hypothesis of stationarity, the distribution of $\hat{ D}_T^*$ is 'close' to the distribution of the random variable $\hat D_T$. On the other hand, under the alternative, $\hat{ D}_T^*$ corresponds to a stationary process with spectral density $g$, which by Theorems \ref{hauptsatz} and \ref{hauptsatz2} implies $\hat{ D}_T^*=O_P(N^{-\gamma})$ for any $\gamma\in(0,1/2)$ while $\hat{ D}_T$ becomes larger than some positive constant, which implies consistency of a bootstrap test. These heuristic arguments are made rigorous by the following statements.
\begin{satz}\label{theorembootstrap}
Let the assumptions of Theorem \ref{hauptsatz} and Assumption \ref{annahmenbootstrap} be fulfilled.  Furthermore, suppose that the following conditions on the growth rate of $p=p(T)$, the estimates $\boldsymbol{\hat{a}}_{j,p}$ and the true AR parameters $\boldsymbol{{a}}_{j}$, $\boldsymbol{{a}}_{j,p}$ defined in  \eqref{Ar1} and \eqref{arpfitt}  are satisfied:
\begin{itemize}
\item[i)] There exist sequences $p_{min}(T)$ and $p_{max}(T)$ such that the order $p$ of the fitted autoregressive process satisfies $p=p(T) \in [p_{min}(T),p_{max}(T)]$ with $p_{max}(T)\geq p_{min}(T)\rightarrow \infty$ and
$
p^{3}_{max}(T)\sqrt{\log(T)/T}=O(1). 
$
\item[ii)] The estimators $\boldsymbol{\hat{a}}_{j,p}$ satisfy
$
\max_{1\leq j\leq p(T)}||\boldsymbol{\hat{a}}_{j,p}-\boldsymbol{a}_{j,p}||_{\infty}=O_P(\sqrt{\log(T)/T}). \label{ARparameterRate}
$
uniformly with respect to $p\in[p_{min}(T),p_{max}(T)]$.
\item[iii)] The matrices $\hat\Sigma_p$ and $\Sigma$ satisfy $\|\boldsymbol{\hat{\Sigma}}_{p}-\boldsymbol{\Sigma}\|_{\infty}\stackrel{P}{\longrightarrow}0$.
\end{itemize}
Then, as $T \rightarrow \infty$, we have 
$
\{\sqrt{N}\hat{\boldsymbol{ D}}_T^*(v,\omega)\}_{(v,\omega)\in[0,1]^2}\Rightarrow\{\boldsymbol{G}(v,\omega)\}_{(v,\omega)\in[0,1]^2}
$
conditionally on $X_{1,T},...,X_{T,T}$,
where $\boldsymbol{G}(v,\omega)$ is the limiting Gaussian process introduced in Theorem \ref{hauptsatz}.
\end{satz}

\begin{satz}\label{theorembootstrap2}
Let the assumptions of Theorem \ref{theorembootstrap}  be fulfilled where (i)  is replaced by
$
p^2_{max}(T)\sqrt{\log(T)/T } N^{1+\epsilon}=o(1)$ { and }  $ {N^{1+\epsilon}}/{p_{\text{min}}(T)}=o(1).
$ 
Then there exists a sequence of $\mathbb{C}^{d\times d}$ valued random processes $\{\hat{\boldsymbol D}_{T,a}^*(v,\omega)\}_{(v,\omega)\in[0,1]}$ such that the following statements hold:
\begin{itemize}
\item[a)] If $K=0$ then for any $T\in\mathbb{N}$: $\sup_{(v,\omega) \in [0,1]^2}||\hat{\boldsymbol D}_{T}(v,\omega)||_{\infty} \stackrel{D}= \sup_{(v,\omega) \in [0,1]^2}||\hat{\boldsymbol D}_{T,a}^*(v,\omega)||_{\infty}$.
\item[b)] If $K\geq 0$, then (as $T\rightarrow\infty$)
\begin{equation}
\frac{ \sup_{(v,\omega) \in [0,1]^2} || \hat{\boldsymbol D}_T^*(v,\omega)||_{\infty}-\sup_{(v,\omega) \in [0,1]^2} || \hat{\boldsymbol D}_{T,a}^*(v,\omega)||_{\infty}}{\E \Big( \sup_{(v,\omega) \in [0,1]^2} || \hat{\boldsymbol D}_{T,a}^*(v,\omega)^2||_{\infty} \Big)^{-1/2}}=o_P(1).
\end{equation}
\item[c)] For any $0<\gamma<\frac{1}{2}$: $N^{\gamma}  \sup_{(v,\omega) \in [0,1]^2} ||\hat{ \boldsymbol D}_{T,a}^*(v,\omega) ||_{\infty} =o_P(1)$.
\end{itemize}
\end{satz}

Note that all assumptions of Theorems \ref{theorembootstrap} and \ref{theorembootstrap2} are rather standard in the framework of AR($\infty$) bootstrap [see for example \cite{bergpappolitis2010} or \cite{kreisspappol2011} among others] and that (ii) is proved in \cite{hannan} for the least squares and the Yule-Walker estimators.
In order to obtain a consistent level $\alpha$ test for the null hypothesis \eqref{null}, we now proceed as follows:\medskip\\
\textbf{Step I (testing for structural breaks)}\label{algorithmtest} We 
calculate the test statistic $\hat D_T$ defined in (\ref{supstat}) and  fit an AR($p$) model to the observed data
 $\{\boldsymbol{X}_{1,T},...,\boldsymbol{X}_{T,T}\}$. 
In the next step Algorithm \ref{algo1} with the corresponding estimates $(\boldsymbol{\hat{a}}_{1,p},...,\boldsymbol{\hat{a}}_{p,p},\boldsymbol{\hat\Sigma}_p)$ 
is used to calculate the bootstrap sample $\{\hat D_{T,1}^*,...,\hat D_{T,B}^*\}$ and the null hypothesis of no break points is rejected, whenever
\vspace{-.5cm}
\be \label{test} \hat D_T >\hat D^*_{T,(\lfloor(1-\alpha)\rfloor B)},
\ee
where $\hat D^*_{T,(\lfloor(1-\alpha)\rfloor B)}$ is an 
 estimate for the $(1-\alpha)$-quantile of the distribution of $\hat D_T$  (here $\hat D^*_{T,(1)},...,\hat D^*_{T,(B)}$ denotes the ordered bootstrap sample). \medskip

Theorem \ref{theorembootstrap} and the continuous mapping theorem imply that the test constructed in (\ref{test}) has asymptotic level $\alpha$ under the assumption $N/T\rightarrow 1/c$. On the other hand in the case $N/T\rightarrow 0$
it can be show by application of Theorem \ref{theorembootstrap2} that the Mallows distance between the random variables 
\vspace{-1cm}
\bea 
&& \sup_{(v,\omega) \in [0,1]^2}||\hat{\boldsymbol D}_T(v,\omega)||_{\infty}\Big/\sqrt{\E  \sup_{(v,\omega) \in [0,1]^2} || \hat{\boldsymbol D}_{T}(v,\omega)^2||_{\infty} }, \\
&&
\sup_{(v,\omega) \in [0,1]^2}||\hat{\boldsymbol D}_T^*(v,\omega)||_{\infty}\Big/\sqrt{\E  \sup_{(v,\omega) \in [0,1]^2} || \hat{\boldsymbol D}_{T}^*(v,\omega)^2||_{\infty} } 
\eea
converges to zero in probability. Therefore similar arguments as in \cite{paparoditis2010} indicate that the bootstrap test has asymptotic level $\alpha$ if $N/T\rightarrow 0$. Moreover, the test is consistent under the alternative in both cases, since it follows from Theorem \ref{theorembootstrap} and \ref{theorembootstrap2} that the bootstrap statistic $\sup_{(v,\omega) \in [0,1]^2} ||\hat{ \boldsymbol D}_T^*(v,\omega)||_{\infty}$ converges to zero in probability, while $\sup_{(v,\omega )\in [0,1]^2} ||  \hat{\boldsymbol D}_T(v,\omega) ||_{\infty}$ becomes larger than some positive constant under the alternative  due to Theorem \ref{hauptsatz}b) and \ref{hauptsatz2}b). The finite sample properties of this test will be investigated in Section \ref{sec6}.
\begin{rem} 
{\rm
If the statistic in \eqref{supstat} is replaced by $\sup\limits_{v\in[0,1]}||\hat{\boldsymbol{D}}_T(v,1)||_{\infty}$, a similar analysis can be performed and we obtain as a special case a test for the null hypothesis that the covariance matrix of the process is constant over time as considered in \cite{aue2009}. Moreover, several other interesting hypotheses can be included by choosing an appropriate function $\boldsymbol \phi(u,\lambda): [0,1] \times [0,\pi] \rightarrow \C^{d \times d}$ and considering functionals of the form
\bea
\frac{1}{e}\left(\int_{0}^{\omega\pi}\int_{v}^{v+e}\boldsymbol \phi(u,\lambda)\boldsymbol{f}(u,\lambda)dud\lambda-\int_{0}^{\omega\pi}\int_{v-e}^{v}\boldsymbol \phi(u,\lambda)\boldsymbol{f}(u,\lambda)dud\lambda \right)
\eea
instead of $\boldsymbol{D}(v,\omega)$. For example, the choice $\boldsymbol \phi(u,\lambda)=\boldsymbol I_{d }\exp(-i\lambda k) $ yields a test for constancy of the  covariance function at a specific lag $k \in \N$. A hypothesis of this type is of interest if the statistician knows in advance that only certain lags influence the dependency structure of the underlying process.
}\end{rem}
\begin{rem}
{\rm
Statistics with a similar structure as $\hat{\boldsymbol D}_T(v,\omega)$ have been studied under different conditions   by several authors in the past  and some references can be found in \cite{giraitis1990}. These authors consider one dimensional time series and the periodogram $I_{(k_1,...,k_2)}(\lambda)$ of the data $X_{k_1},X_{k_1+1},...,X_{k_2}$. In order to test for a change point they propose to compare the estimators $\hat F_k(\omega)=\int_0^{\omega \pi} I_{(1,...,k)}(\lambda)$ and $\hat F_{T-k}^*(\omega)=\int_0^{\pi \omega}I_{(k+1,...,T)}(\lambda)d\lambda$ by  calculating 
$
\sup_{k \in \{1,...,T\}}\sup_{\omega \in [0,1]}|\hat F_k(\omega)-\hat F_{T-k}^*(\omega)|.
$
Note that this procedure can detect at most one break point and that these authors show that an appropriately standardized test statistic converges to some non degenerate limit whose quantiles are, however, unknown. These results are of great importance to understand the theoretical properties of such statistics. On the other hand, the construction of a computable level $\alpha$ test based on statistics of this type is an open and very challenging problem.}
\end{rem}

\begin{rem}\label{gaussian}
{\rm
We emphasize that the assumption of Gaussianity in this (and the following) section is merely imposed to simplify technical arguments in the proofs. In particular it is straightforward (but cumbersome) to extend the results to a more general class of linear processes. In fact, the proof of Theorem \ref{hauptsatz} can be modified to address for non Gaussian innovations. The only difference appears in the bootstrap procedure because Gaussian  innovations are not an appropriate choice of the bootstrap replicates $\boldsymbol{Z}_t^*$ in  Algorithm \ref{algo1}. In general one has to employ replicates $\boldsymbol{Z}_t^*$ which mimic the fourth cumulant of the true underlying innovation process [see \cite{krepap2012} for more details].
}
\end{rem}

\section{Detecting the number and location of break points} \label{sec4}
\def\theequation{4.\arabic{equation}}
\setcounter{equation}{0}
If structural breaks have been detected by the test \eqref{test} it is of futher interest to estimate the number and location of possible break points and to identify the components responsible for these changes in the regime. In the following discussion we will develop a procedure which consists of two further steps and detects simultaneously the number, location and corresponding components of multiple structural breaks. In the second step we estimate  (shrinking if $N/T\rightarrow 0$) sets, which may contain potential break points. Roughly speaking these sets contain all points where the components of the spectral density estimate indicate a structural break. In the third step these sets are reduced to identify the relevant structural breaks and corresponding components which are responsible for these breaks. For this purpose we recall the definition 
\eqref{empprocess},
choose some constant $0 < \gamma <1/2$ (a recommendation for this choice will be given in Section \ref{sec5}) and proceed as follows.
\medskip

\textbf{Step II (identification of sets containing break points)}
We consider a point $v \in \{\frac{N}{T},\frac{N+1}{T},...,\frac{T-N}{T}\}$ as a candidate for a structural break in the component $(a,b)$ if the inequality
\vspace{-.3cm}
\begin{equation}
N^{\gamma}\sup\limits_{\omega \in [0,1]}|\hat{[\boldsymbol{D}}_T(v,\omega)]_{a,b}|>\boldsymbol{\epsilon}_{T,a,b}(v)\label{QQ2}
\end{equation}
holds, where  $\boldsymbol{\epsilon}_{T,a,b}(v)$ is a threshold satisfying $\liminf\limits_{T\rightarrow\infty} \boldsymbol{\epsilon}_{T,a,b}(v)\geq C>0$ for some constant $C$ and $\boldsymbol{\epsilon}_{T,a,b}(v)=o(N^{\gamma})$ uniformly in $v \in [0,1]$. 
A data driven rule for the choice of the threshold $\boldsymbol{\epsilon}_{T,a,b}(v)$ with good finite sample properties will be given in Section \ref{sec5}. \medskip

The decision rule (\ref{QQ2}) identifies subsets $R_1,...,R_{K_T}\subset\{N/T,....,1-N/T\}$ where possible break points in the components of the spectral density matrix may occur.
The following example illustrates the second step of the procedure.
\begin{example}{\rm \label{example1}
We consider the bivariate model
\vspace{-.5cm}
\begin{align} \label{beispielbild}
&\boldsymbol{X}_{t,T}=\sum_{j=1}^41_{(\frac{j-1}{4}T,\frac{j}{4}T]}(t)\Theta_j\boldsymbol{Z}_t,
\end{align}
where the matrices $\Theta_1,\Theta_2,\Theta_3,\Theta_4$ are defined by  \vspace{-1cm}
\begin{singlespace}
\begin{align*}
\Theta_1:=
\begin{pmatrix}
1 & 0\\
0 & 1
\end{pmatrix}\quad
\Theta_2:=
\begin{pmatrix}
2 & 0\\
0 & 1
\end{pmatrix}\quad
\Theta_3:=
\begin{pmatrix}
2 & 0\\
0 & 2
\end{pmatrix}\quad
\Theta_4:=
\begin{pmatrix}
\sqrt{2} & \sqrt{2}\\
0 & 2
\end{pmatrix}
\end{align*}
\end{singlespace}
and $\{\boldsymbol{Z}_t\}_{t\in\mathbb{Z}}$ is a two dimensional Gaussian white noise process.
The spectral density $\boldsymbol{f}$ of a bivariate time series model \eqref{beispielbild} exhibits $3$ break points, where the first change only involves the first component $[\boldsymbol{f}]_{1,1}$, the second  only concerns the component $[\boldsymbol{f}]_{2,2}$ and the third break point leaves the compontens $[\boldsymbol{f}]_{1,1}$ and $[\boldsymbol{f}]_{2,2}$ unchanged but appears in the 
cross spectrum  $[\boldsymbol{f}]_{1,2}$ and $[\boldsymbol{f}]_{2,1}$.
Figure \ref{illustration1} contains a plot of a typical set of data of length $T=2048$ generated by model \eqref{beispielbild}, 
{where the dashed vertical lines indicate the true break points in the univariate time series. Note that the third break point only 
corresponds to a change in the dependency structure of the two univariate data sets.} Figure \ref{illustration2} shows the four plots of the functions $v \mapsto N^{\gamma}\sup\limits_{\omega \in [0,1]}|[\hat{D}_T(v,\omega)]_{a,b}|$, $a,b \in \{1,2\}$, (solid lines), where $N=256$, $\gamma=0.49$. In each component we added a plot of the threshold level $v \mapsto \boldsymbol{\epsilon}_{T,a,b}(v)$,  $a,b \in \{1,2\}$ (broken lines), which will be formally defined in equation (\ref{epsilonchoice}) of the following section. It is evident that for each component $(a,b)$ the test statistic exceeds the level $\boldsymbol{\epsilon}_{T,a,b}(v)$ in a neighborhood of the break point. For the simulated scenario we obtain the sets $R_1=\{\frac{289}{2048},...,\frac{606}{2048}\}$, $R_2=\{\frac{802}{2048},...,\frac{1140}{2048}\}$, $R_3=\{\frac{1427}{2048},...,\frac{1715}{2048}\}$ of potential break points. Note that the local maximum of the function $N^{\gamma}\sup\limits_{\omega \in [0,1]}|[\hat{\boldsymbol{D}}_T(v,\omega)]_{a,b}|$ is rather close to the true change point. }
\begin{figure}[]
\begin{center}
\includegraphics[width=150mm,height=70mm]{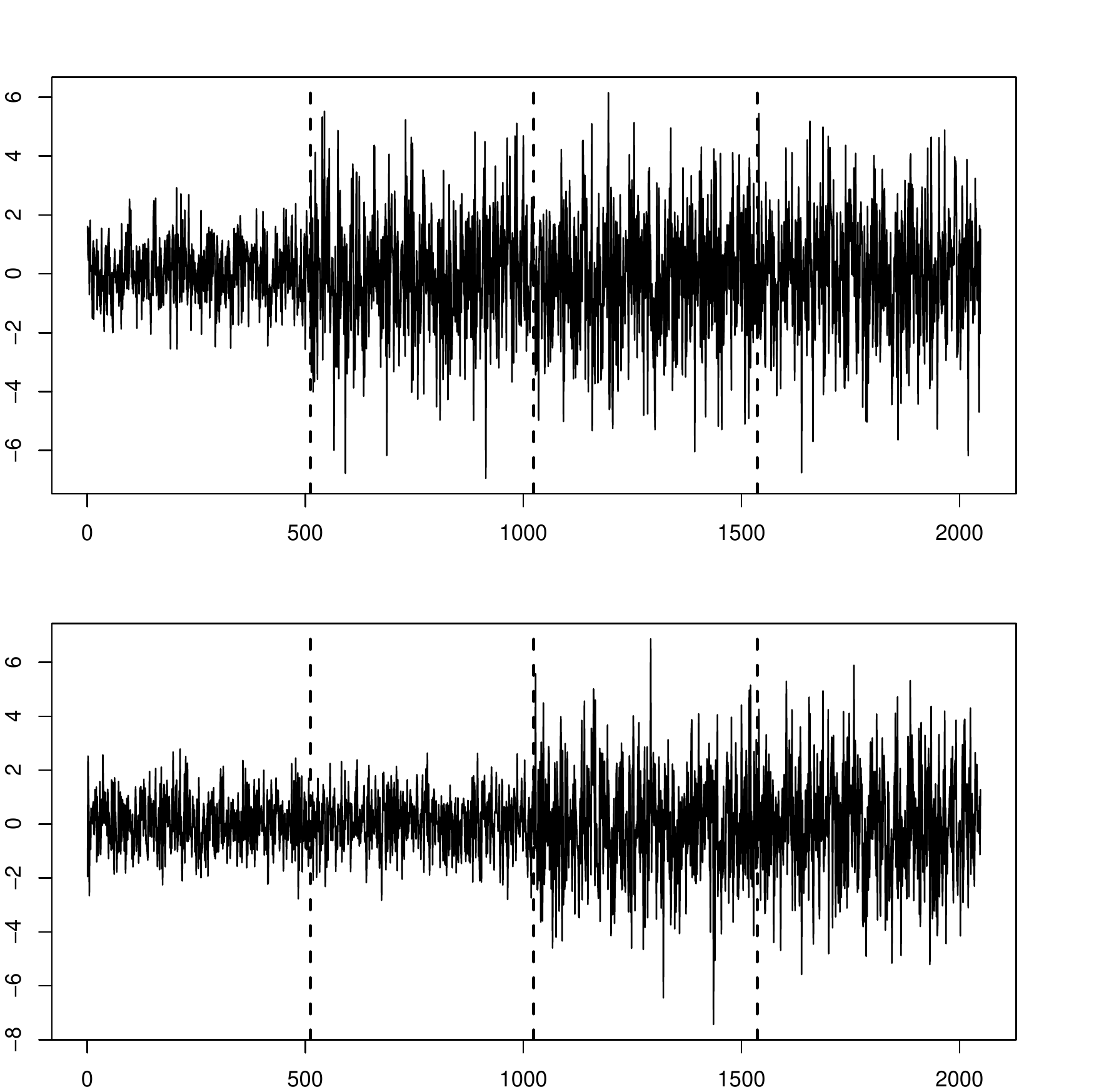}
\caption{\textit{ Simulated data from the model \eqref{beispielbild}, where $T=2048$. {The vertical dashed lines denote the true break points in the univariate time series.}}} 
\label{illustration1}
\end{center}
\end{figure}
\begin{figure}[]
\begin{center}
\includegraphics[width=170mm,height=80mm]{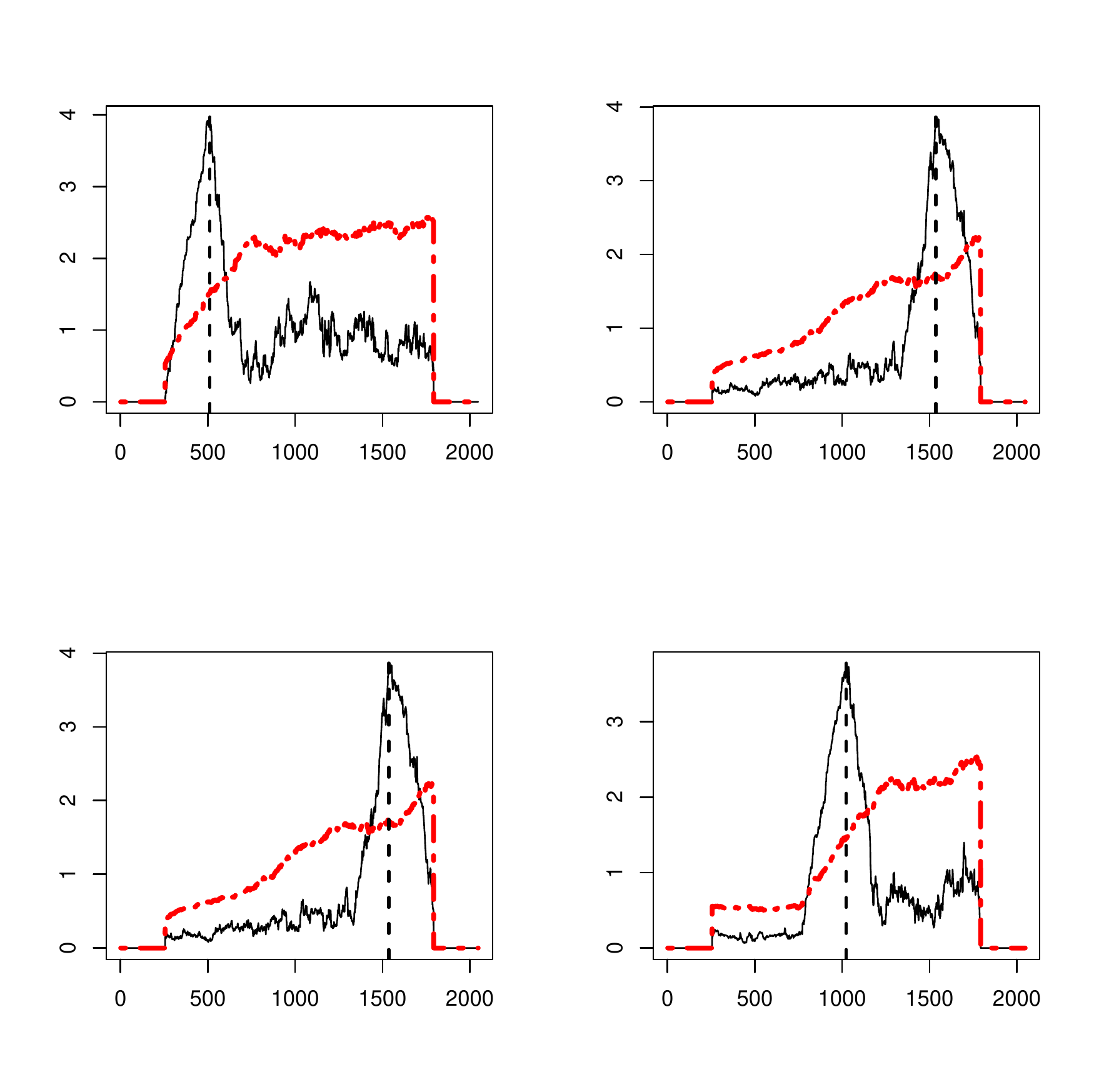}
\caption{\textit{Plots of the functions $v \mapsto N^\gamma \sup\limits_{\omega \in [0,1]}\big|[\hat{D}_T(v,\omega)]_{a,b}\big|$ (solid lines) and $v \mapsto \boldsymbol{\epsilon}_{T,a,b}(v)$ $(a,b=1,2)$ (broken lines) with vertical dashed lines at the true break points.}} 
\label{illustration2}
\end{center}
\end{figure}
\end{example}
The following result shows that for an increasing sample size the subsets $R_1,...,R_{K_T}$ are contained in neighborhoods of radius of $\frac{N}{T}$ of the 'true' break points, that is
\vspace{-.7cm}
\bea
\bigcup_{j=1}^{K_T}R_j ~\subset ~~\bigcup_{a,b \in \{1,...,d\}} \mathcal{I}_{T,a,b}(b_1,...,b_K),
\eea
 \vspace{-.5cm}
where
\begin{equation}\label{tauab}
\mathcal{I}_{T,a,b}(b_1,...,b_K):=\bigcup\limits_{\substack{j=1 \\ \sup_\omega |[D(b_j,\omega)]_{a,b}| >0 }}^K\Big\{\frac{\lfloor b_jT\rfloor-N}{T},...,\frac{\lfloor b_jT\rfloor+N}{T}\Big\}.
\end{equation}
\begin{satz}\label{detection}
Assume that condition \eqref{summ1} is satisfied and that the sequence $N$ satisfies one of the following conditions:
\begin{itemize}
\item[(1)] There exists some $\epsilon>0$ such that $N\rightarrow\infty$, $T^{\epsilon}/N\rightarrow 0$ and $N/T\rightarrow0$.
\item[(2)] There exists some constant $c\geq 2/\min\limits_{i=1,...,K-1}|b_i-b_{i+1}|$ such that $N/T \rightarrow 1/c$.
\end{itemize}
Furthermore for all $a,b \in \{1,...,d\}$, $v\in[0,1]$ and $\gamma\in(0,\frac{1}{2})$ let $(\boldsymbol{\epsilon}_{T,a,b}(v))_{T\in\N}$ denote a sequence satisfying $\boldsymbol{\epsilon}_{T,a,b}(v)=o(N^{\gamma})$ and $\liminf\limits_{T\rightarrow\infty} \inf\limits_{v \in [0,1]} \boldsymbol{\epsilon}_{T,a,b}(v)\geq C>0$ for some constant $C$. Then the detection rule (\ref{QQ2}) is accurate in the following sense:
\begin{itemize}
\item[a)] The probability that the decision rule (\ref{QQ2}) indicates a structural break at a rescaled time point, which has a distance of at least $\frac{N}{T}$ from each  of the break points $b_1,...,b_K$, vanishes asymptotically, that is
\vspace{-.5cm}
\begin{equation}\label{falselydetecting}
\mathbb{P}\Big(\bigcup_{a,b \in \{1,...,d\}}\bigcup\limits_{v \in \bar{\mathcal{I}}_{T,a,b}(b_1,...,b_K)}\Big\{N^{\gamma}\sup\limits_{\omega \in [0,1]}|[\hat{\boldsymbol{D}}_T(v,\omega)]_{a,b}|>\boldsymbol{\epsilon}_{T,a,b}(v)\Big\}\Big)\xrightarrow{\text{ }T \rightarrow \infty\text{ }} 0,
\end{equation} 
where $\bar{\mathcal{I}}_{T,a,b}(b_1,...,b_K)=\{\frac{N}{T},\frac{N+1}{T},...,\frac{T-N}{T} \}\setminus\mathcal{I}_{T,a,b}(b_1,...,b_K)$.
\item[b)] The probability that the procedure detects all structural breaks  converges to one, that is
\vspace{-.5cm}
\begin{equation}
\mathbb{P}\Big(\bigcap\limits_{v\in\{b_1,...,b_{K}\}}\bigcap\limits_{(a,b)\in B(v)}\Big\{N^{\gamma}\sup\limits_{\omega \in [0,1]}|[\hat{\boldsymbol{D}}_T(v,\omega)]_{a,b}|>\boldsymbol{\epsilon}_{T,a,b}(v)\Big\}\Big)\xrightarrow{\text{ }T \rightarrow \infty\text{ }} 1,
\end{equation}
where $B(v):=\{(a,b) \in \{1,...,d\}^2| \sup_{\omega \in [0,1]} |[\boldsymbol{D}(v,\omega)]_{a,b}| > 0 \}.$
\end{itemize}
\end{satz}

Recall that we use the decision rule \eqref{QQ2} to identify sets of possible break points. By Theorem \ref{detection} it follows with probability converging to $1$ that these sets are contained in the (shrinking if $N/T\rightarrow 0$) neighborhoods $\{\frac{\lfloor b_jT\rfloor-N}{T},...,\frac{\lfloor b_jT\rfloor+N}{T}\}$ of those points $b_j$, for which there is at least one change in one of the components of the spectral density matrix. As demonstrated in Example \ref{example1} the inequality \eqref{QQ2} is usually satisfied in a neighborhood of the break point. Consequently, for a finite sample size, the number of possible change points detected by \eqref{QQ2} is usually much larger than the true number $K \in \N$. To address this issue these sets are further reduced in the third step of the procedure. As a final result we obtain estimators of the number and the location of structural breaks. The basic idea is very simple: For a certain set $R_j$ of points in $\{N/T,...,1-N/T\}$ satisfying \eqref{QQ2} for at least one pair $(a,b)\in\{1,...,d\}^2$ we identify a point $\tilde{b}\in R_j$ for which the local deviation from stationarity is maximal and then remove all points of the interval $[\tilde{b}-\frac{N}{T},\tilde{b}+\frac{N}{T}]$ from the set $R_j$. The detailed description of this idea, which has to be iterated, is described below. \medskip

\textbf{Step III (localization of structural breaks)}
\begin{itemize}
\vspace{-.1cm}
{
\item[(1)] Let $\hat K$ denote the number of elements $v\in\{\frac{N}{T},\frac{N+1}{T},...,\frac{T-N}{T}\}$ which fulfill \eqref{QQ2} for at least one component $(a,b)\in\{1,...,d\}^2$, denote the corresponding elements by $\hat b_1,...,\hat b_{\hat K}$ and define the sets $\hat B_P:=\{\hat b_1,...,\hat b_{\hat K}\}$ and $\hat B_D=\emptyset$ of potential break points and detected break points, respectively.
\vspace{-.3cm}
\item[(2)] If the set $\hat B_P$ is not empty add the element $\tilde b \in \hat B_P$ to the set $\hat B_D$ for which the statistic
$
\sup_{(a,b) \in \{1,...,d\}^2} ( \sup\limits_{\omega \in [0,1]}N^{\gamma}|[\hat{\boldsymbol{D}}_T(\tilde b,\omega)]_{a,b})
$
is maximal and replace the set $\hat B_P$ by $\hat B_P \backslash [\tilde b-\frac{N}{T},\tilde b+\frac{N}{T}]$.
\vspace{-.3cm}
\item[(3)] Repeat step (2) until $B_P=\emptyset$ and 
redefine $\hat K=|\hat B_D|$ and $(\hat b_1,...,\hat b_{\hat K})$ as the elements of $\hat B_D$ such that $\hat b_i < \hat b_{i+1}$ for $i=1,...,\hat K-1$.}
\end{itemize}

Note that in step  $(2)$ of the above procedure we pick the element from $\hat B_P$ for which a change in the spectral density matrix seems to be most noticable. If $\hat K$ break points $\hat b_1,...,\hat b_K$ have been found, it is (as mentioned above) of further interest to determine the nature of each detected structural break, i.e. to identify the components of the spectral density matrix, which exhibit a change at the specific points in time. Such a refined analysis can be undertaken by considering a change in the component $(a,b)$ of the spectral density matrix $\boldsymbol{f} \in \C^{d\times d}$  at the point $\hat b_i$ as relevant if the inequality
$
N^{\gamma}\sup\limits_{\omega \in [0,1]}\big|[\hat{\boldsymbol{D}}_T(\hat b_i,\omega)]_{a,b}\big|>\boldsymbol{\epsilon}_{T,a,b}(\hat b_i)
$
holds. 
\begin{example}{\rm
In order to illustrate this step we continue the discussion of Example \ref{example1} where the sets $R_1$, $R_2$ and $R_3$ of potential break points were identified. Therefore at the beginning of Step III the set  of possible break points is given by $\hat B_P=R_1\cup R_2 \cup R_3$ while the set of determined break points satisfies $\hat B_D=\emptyset$. In the first step we add the element $\tilde b=\frac{511}{2048}$ to $\hat B_D$ and then remove all elements, which do not have a distance of at least $N/T=\frac{256}{2048}$ from $\tilde b=\frac{511}{2048}$, from the set $\hat B_P$. This leaves the set $\hat B_P= R_2\cup R_3$. In the next iteration we add the element $\tilde b=\frac{1537}{2048}$ to $\hat B_D$ and reduce the set $\hat B_P$ to $\hat B_P= R_2$. In the last step we add the point $\tilde b=\frac{1026}{2048}$ and obtain $\hat B_D=\{\frac{511}{2048},\frac{1026}{2048},\frac{1537}{2048}\}$ and reduce the set $\hat B_P$ accordingly. At this point the procedure terminates because $\hat B_P=\emptyset$ and yields  $\hat K=3$ as an estimator for the number of break points. The locations of the break points are estimated as $(\hat b_1,\hat b_2,\hat b_3)=(\frac{511}{2048},\frac{1026}{2048},\frac{1537}{2048})$. Thus in this example all break points are detected and the respective components of the spectral density matrix $\boldsymbol{f}$ responsible for these changes are identified as well. } 
\end{example}
Theorem \ref{detection} implies that $\hat K$ converges to $K$ in probability and that the same statement also holds for the estimated break points $\hat b_1,...,\hat b_{\hat K}$, that is
\vspace{-.5cm}
\be\label{consis}
\hat b_j\stackrel{P}{\longrightarrow}b_j\quad (j=1,...,K),
\ee
where $b_1,...,b_K$ denotes the true locations of the break points. In the case (1) of shrinking relative size of the blocklength $N/T \to 0$, the property \eqref{consis} simply follows from the fact that the diameter of the sets $R_j$ converges to zero. For sequences $N$ which satisfy condition (2), consistency follows from the property 
\vspace{-.5cm}
\be\label{propertydetection}
N^{\gamma}\sup\limits_{(v,\omega)\in[0,1]^2}\|\hat {\boldsymbol{D}}_T(v,\omega)-\boldsymbol{D}_{N,T}(v,\omega)\|_{\infty}=o_P(1),
\ee
where 
$\boldsymbol{D}_{N,T}(v,\omega):=\frac{T}{N} (\int_{0}^{\omega\pi}\int_{v}^{v+N/T}\boldsymbol{f}(u,\lambda)dud\lambda-\int_{0}^{\omega\pi}\int_{v-N/T}^v\boldsymbol{f}(u,\lambda)dud\lambda ).$
This  can be shown by similar arguments as given in the proof of Theorem \ref{hauptsatz} and the fact that the deterministic sequence $\bar D_{N,T}(v):=\sup\limits_{\omega\in[0,1]}\|\boldsymbol{D}_{N,T}(v,\omega)\|$  attains local maxima at the true break points $b_j$ as $T \rightarrow \infty$.
In addition, Theorem \ref{detection} also shows that the refined analysis described in the previous paragraph consistently identifies the components where the structural breaks are present. This means that the procedure is able to detect the 'true' number, locations and components responsible for the breaks with a probability converging to one as the sample size increases. We will give data driven rules for the choice of the regularization parameters in Section \ref{sec5} and investigate the finite sample properties of the  procedure in Section \ref{sec6}.

\section{Data driven regularization} \label{sec5}
\def\theequation{5.\arabic{equation}}
\setcounter{equation}{0}

The proposed test for structural breaks and the corresponding localization method  depend on the choice of several regularization parameters. 
In   this section we will present guidelines for a data driven choice of these parameters which from a theoretical point of view satisfy the assumptions for the asymptotic theory discussed in the previous sections and from a practical point of view, yield good results for finite samples (see the examples of the following section).\\
We begin with 
 the choice of $\gamma$, where  one faces a tradeoff in practical applications. Roughly speaking large values of $\gamma\in(0,\frac{1}{2})$ tend to increase the number $\hat K$ of estimated break points. This property is of practical importance for example if it is desirable to avoid an underestimation of the true number of break points. In our simulation study we follow this path (i.e. we want to be sure to detect all underlying change points) and choose $\gamma=0.49$. This value leads to satisfactory results in all investigated scenarios (see Section \ref{sec6} for more details) and is therefore recommended. \\
For the  choice of the threshold $\boldsymbol{\epsilon}_{T,a,b}(v)$ in \eqref{QQ2}  we note that it follows by similar calculations as given in the proof of
Theorem \ref{hauptsatz}   that the  variance of $\hat{[\boldsymbol{D}}_T(v,\omega)]_{a,b}$ satisfies
\vspace{-.5cm}
\begin{align*}
\lim\limits_{\substack{T\rightarrow\infty\\N/T\rightarrow c}}\Var(&\sqrt{N}[\hat{\boldsymbol{D}}_T(v,\omega)]_{a,b})  = \frac{1}{2\pi} \int_0^{\pi \omega}  \Big( [\boldsymbol{f}(v-\frac{1}{2c},\lambda)]_{aa}[\boldsymbol{f}(v-\frac{1}{2c},\lambda)]_{bb}\\
+&[\boldsymbol{f}(v+\frac{1}{2c},\lambda)]_{aa}[\boldsymbol{f}(v+\frac{1}{2c},\lambda)]_{bb}+ |[\boldsymbol{f}(v-\frac{1}{2c},\lambda)]_{ab}|^2+|[\boldsymbol{f}(v+\frac{1}{2c},\lambda)]_{ab}|^2 \Big) d\lambda.
\end{align*}
Consequently, for fixed $v,\omega$, the asymptotic local variance can be easily  estimated on a block of length $2N$ by 
$
M_{T,a,b} (v,\omega) =
{\frac{1}{N}\sum_{k=1}^{  \lfloor \omega N\rfloor    }[\boldsymbol{I}_{2N}(v,\lambda_{k,2N})]_{aa}[\boldsymbol{I}_{2N}(v,\lambda_{k,2N})]_{bb}}, 
$
where $\lambda_{k,2N}:=2\pi k/2N$. 
Similar arguments as given in the proof of Theorem \ref{hauptsatz} yield that $N^{\gamma}[\hat{\boldsymbol{D}}_T(v,\omega)]_{a,b}$ and $N^{\gamma}[\hat{\boldsymbol{D}}_T(v^\prime,\omega^\prime)]_{a,b}$ are asymptotically independent whenever $v\neq v^\prime$. Motivated by hard thresholding [see \cite{thresholding2}] 
we recommend
\vspace{-.5cm}
\be
\boldsymbol{\epsilon}_{T,a,b}(v)&=& \sqrt{2M_{T,a,b}(v,1)\log\Big(\frac{d(d+1)T}{2N}\Big)} \label{epsilonchoice} ,
\ee
where $d$ denotes the dimension of the time series under consideration.
 Note that under the condition $\inf\limits_{u\in[0,1]}\min\limits_{i\in\{1,...,d\}}[\boldsymbol{f}(u,\lambda)]_{ii}>0$ it follows that $\liminf\limits_{T\rightarrow\infty}\boldsymbol{\epsilon}_{T,a,b}(v)\geq C>0$ for all $v\in[0,1]$, $a,b\in\{1,...,d\}$ with probability converging to 1, and that the threshold sequence \eqref{epsilonchoice} thus satisfies the conditions of Theorem \ref{detection}.\\
The choice of the window length $N$ has to reflect two contradicting objectives. On the one hand, $N$ should be chosen rather large, if there exists no change point or if there are only a few structural breaks with long stationary segments in between. On the other hand, $N$ should be chosen rather small if there exist many break points with small distances. This basic idea is incorporated in the following rule.
\begin{algo}\label{algorithmchoicen}{ \rm {\bf (choice of the window length $N$)}
We consider a set of even integers, say  $V_T:=\{N_1,...,N_{n}\}$ 
satisfying $\sqrt{T}\leq N_1<N_2<\hdots <N_{n}\leq T^{5/6}$ and determine for
 each $N\in V_{T}$   the number $\hat K_T(N)$ of break points detected by the algorithm of Step III.
 We define $i^*:=\sup\{i\in\{2,...,n(T)\}|\hat K_T(N_{i-1})\leq\hat K_T(N_{i})\}$ (here $\sup \emptyset = - \infty$),
 $N^*=N_{i^*}$ if $i^*\leq n(T)$ and $N^*=N_{n(T)}$ if $i^*=-\infty$
and use $N=2N^*$ for the test of structural breaks and $N=N^*$ for the estimation of $K$ and the localization of the break points.
}
\end{algo}

Algorithm \ref{algorithmchoicen} employs the detection method described in Section 4 for every $N\in V_T$ and selects $N^*$ as the largest $N\in V_T$ for which there is no additional break point detected for the next smaller $N\in V_T$. Note also that the recommended choice of $N$ is different for  step I and steps II, III of the procedure. This recommendation is motivated by the following observation. If the distance between two consecutive break points is $m$, then one should use $N=m$ in step I, while the 'best' window length would be $N=m/2$ in step II [because  we  assume
 that the distance between two consecutive change points is at least $2/c$]. In the numerical investigations we additionally restrict ourselves to window lengths which equal a power of two, that is we consider $n=\lfloor\log_2(T^{5/6})\rfloor-\lceil\log_2(\sqrt{T})\rceil+1$ and $N_i=2^{\lceil\log_2(\sqrt{T})\rceil-1+i}$ for $i=1,...,n$. This choice allows for an application  of the FFT algorithm in the calculation of the local periodogram which yields a significant reduction in computational time. \\
 Finally we  select  the order $p$ for the autoregressive bootstrap  as the minimizer of the AIC criterion, which is defined by
 \vspace{-.5cm}
\begin{align*}  
\hat p = \text{argmin}_p \frac{2\pi}{T} \sum_{k=1}^{{T}/{2}} \Big( \log(\text{det}[\boldsymbol{f}_{\hat \theta(p)}(\lambda_{k,T})])+\text{tr}[(\boldsymbol{f}_{\hat \theta(p)}(\lambda_{k,T}))^{-1}\boldsymbol I_T(\lambda_{k,T})] \Big) + p/T
\end{align*}
in the context of stationary processes [see \cite{whittle1}]. Here, $\boldsymbol f_{\hat \theta(p)}$ is the spectral density of the fitted stationary AR($p$) process and $\boldsymbol I_T$ is the usual periodogram calculated under the assumption of stationarity with the corresponding Fourier frequencies $\lambda_{k,T}=2\pi k/T$.

\section{Finite sample properties}  \label{sec6}
\def\theequation{6.\arabic{equation}}
\setcounter{equation}{0}
In this section we investigate the finite sample properties of the proposed procedure. We also provide a comparison with competing methods and illustrate the methodology in a data example. The regularizing parameters are chosen as described in the previous section. Throughout this paragraph $\{\boldsymbol{Z_t}\}_{t\in\mathbb{Z}}$ always denotes a standard normal distributed White Noise sequence. \\
We begin with a study of the power of the bootstrap test, where 
$1000$  and $500$ simulation runs are used for estimating the rejection probabilities under the null hypothesis and alternative, respectively, and 
 $300$ bootstrap replications are performed  for the test \eqref{test}.
In order to investigate the approximation of the nominal level we consider the bivariate MA(1) and AR(1) models
\vspace{- 1.25cm}
\begin{singlespace}
\begin{align}
\boldsymbol{X_t}&=\boldsymbol{Z_t}+\begin{pmatrix}
\theta & 0.2 \\
0.2 & \theta
\end{pmatrix}
\boldsymbol{Z_{t-1}}\label{Ma1}\\
\boldsymbol{X_t}&=
\begin{pmatrix}
\phi & 0.2\\
0.2 & \phi
\end{pmatrix}
\boldsymbol{X_{t-1}}+\boldsymbol{Z_t}.\label{Ar1multi}
\end{align}
\end{singlespace}
The simulated rejection frequencies for different values of the parameters $\theta$ and $\phi$ are given in Table \ref{tableH_0} and we observe that the nominal level is underestimated for small sample sizes, but the approximation becomes better with increasing sample size. 
\begin{table}[htbp]\footnotesize
\centering
\begin{tabular}{|c|c|c|c|c|c|c|c|c|}
\hline
\multicolumn{1}{|c|}{} & \multicolumn{4}{c|}{$\text{H}_0$: Model (\ref{Ma1})} & \multicolumn{4}{c|}{$\text{H}_0$: Model (\ref{Ar1multi})} \\ \hline
& \multicolumn{ 2}{c|}{$\theta=0.5$} & \multicolumn{ 2}{c|}{$\theta=-0.5$} & \multicolumn{ 2}{c|}{$\phi=0.5$} & \multicolumn{ 2}{c|}{$\phi=-0.5$} \\ \hline
$T$ & 5\% &10\% & 5\% &10\% & 5\% &10\% & 5\% &10\%  \\ \hline
128 & 0.018 & 0.049 & 0.019 & 0.047 & 0.021 & 0.061 & 0.020 & 0.059 \\ \hline
256 & 0.025 & 0.054 & 0.029 & 0.063 & 0.031 & 0.070 & 0.031 & 0.068 \\ \hline
512 & 0.043 & 0.081 & 0.039 & 0.088 & 0.034 & 0.075 & 0.040 & 0.083 \\ \hline
\end{tabular}
\caption{\textit{Simulated nominal level of the test \eqref{test} for the models \eqref{Ma1} and \eqref{Ar1multi} with different choices of $\theta$, $\phi$ and $T$.}}
\label{tableH_0}
\end{table}
Next we study the power of the test \eqref{test} and compare it with the CUSUM type procedure proposed by \cite{aue2009}. This procedure is specifically designed to test for constancy of the covariance matrix $\boldsymbol{\Cov(X_t,X_t)}$.
%
We consider  three bivariate models
\vspace{-1cm}
\begin{singlespace}
\begin{align}
\boldsymbol{X}_{t,T}=&\sum_{l=0}^K1_{[\lfloor b_lT\rfloor+1,\lfloor b_{l+1}T\rfloor]}(t)
\begin{pmatrix}
\phi_l&0.1\\
0.1&\phi_l
\end{pmatrix}
\boldsymbol{X}_{t-1,T}+\boldsymbol{Z}_t\label{model4b}\\
\boldsymbol{X}_{t,T}=&\sum_{l=0}^K1_{[\lfloor b_lT\rfloor+1,\lfloor b_{l+1}T\rfloor]}(t)
\begin{pmatrix}
\theta_l&0.1\\
0.1&\theta_l
\end{pmatrix}
\boldsymbol{Z}_{t-1}+\boldsymbol{Z}_t\label{model4c}\\
\boldsymbol{X}_{t,T}=&\sum_{l=0}^K1_{[\lfloor b_lT\rfloor+1,\lfloor b_{l+1}T\rfloor]}(t)
\begin{pmatrix}
\sigma_l&0.2\\
0.2&\sigma_l
\end{pmatrix}
\boldsymbol{Z}_t\label{model4}
\end{align}
\end{singlespace}
for different choices of the number $K$ and location $\mathfrak{b}=(b_1,...,b_K)$  of the break points and parameters $\Sigma:=(\sigma_0,...,\sigma_K)$, $\Phi:=(\phi_0,...,\phi_K)$ and $\Theta:=(\theta_0,...,\theta_K)$.
The results are summarized in Table \ref{vglmitAUE1}.

\begin{table}[htbp]\footnotesize
\centering
\begin{tabular}{|c|c|c|c|c|c|c|c|c|}
\hline
\multicolumn{1}{|c|}{} &\multicolumn{2}{|c|}{} & \multicolumn{2}{|c|}{$T=128$} & \multicolumn{2}{|c|}{$T=256$} & \multicolumn{2}{|c|}{$T=512$}\\\hline
& $\mathfrak{b}$ & parameter & \eqref{test} & Aue & \eqref{test} & Aue & \eqref{test} & Aue \\\hline\hline
\multirow{2}{*}{(\ref{model4b})} & $(\frac{1}{4},\frac{2}{3},\frac{3}{4})$ & $\Phi=(0.5,-0.5,0.5,-0.5)$ & 0.072 & 0.046 & 0.306 & 0.059 & 0.812 & 0.090\\
& $(\frac{1}{2})$ & $\Phi=(0.5,-0.5)$ & 0.078 &0.047 & 0.234 & 0.146 & 0.702 & 0.337\\\hline\hline
\multirow{2}{*}{(\ref{model4c})} & $(\frac{1}{4},\frac{2}{3},\frac{3}{4})$ & $\Theta=(1,-1.5,1,-1.5)$ & 0.284 &0.066 & 0.544 & 0.119 & 0.926 & 0.282\\
& $(\frac{1}{2})$ & $\Theta=(1,-1.5)$ & 0.258 & 0.257 & 0.542 & 0.620 & 0.850 & 0.965 \\\hline\hline
\multirow{2}{*}{(\ref{model4c})} & $(\frac{1}{4},\frac{2}{3},\frac{3}{4})$ & $\Sigma=(1,2,1,0.5)$ & 1.000 & 0.213 & 1.000 & 0.989 & 1.000 & 1.000\\
  &$(\frac{1}{2})$ & $\Sigma=(1,2)$ & 0.864 & 0.818  & 0.990 & 1.000 & 1.000 & 1.000 \\
 &$\emptyset$ & $\Sigma=(1)$ & 0.026 & 0.033  & 0.029 & 0.048 & 0.041 & 0.047 \\\hline
\end{tabular}
\caption{\textit{Empirical rejection frequencies of the test \eqref{test} and the CUSUM type procedure of \cite{aue2009}. The nominal level is $\alpha=0.05$ and the models are defined in (\ref{model4b}), (\ref{model4c}) and (\ref{model4}), where different choices of break points $\mathfrak{b}=(b_1,...,b_K)$ and AR-parameters $\Phi=(\phi_0,...,\phi_K)$, MA-parameters $\Theta=(\theta_0,...,\theta_K)$ and standard deviations $\Sigma=(\sigma_0,...,\sigma_K)$ are considered.}}
\label{vglmitAUE1}
\end{table}


We observe from the upper part of Table \ref{vglmitAUE1} that in the AR-model \eqref{model4b} the test proposed in this paper significantly outperforms the CUSUM method of \cite{aue2009} in all cases under consideration. In the MA-model \eqref{model4c} the new test yields a substantially larger power if there exist multiple break points while the power is slightly smaller if there exists only one break point.

\begin{figure}
\begin{center}
\includegraphics[width=0.9\textwidth,height=85mm]{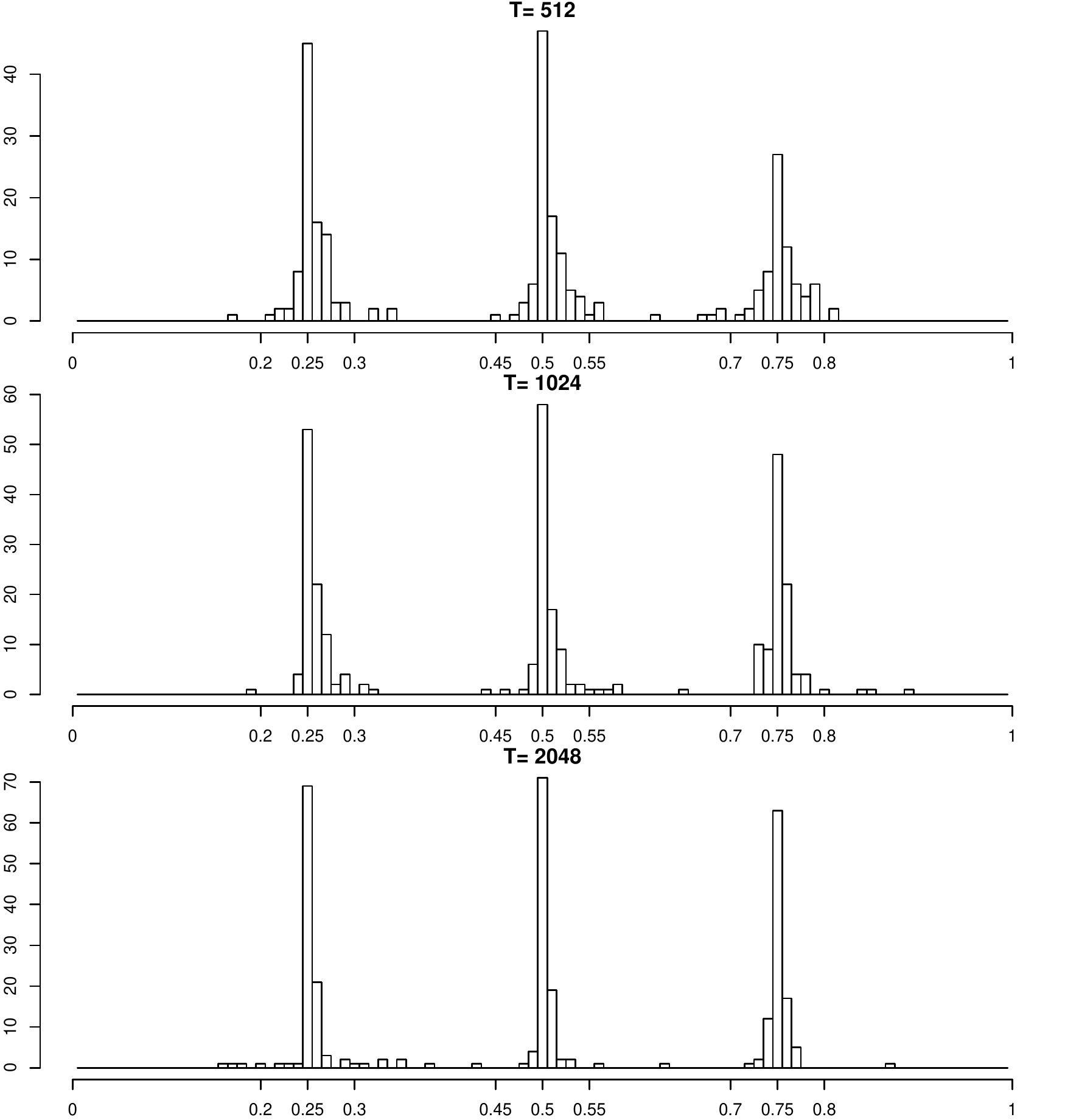}
\caption{\textit{Histograms for the empirical distribution of $\hat{\mathfrak{b}}=(\hat b_1,...,\hat b_{\hat K})$ based on 100 simulation runs of model \eqref{beispielbild} for sample sizes $T\in\{512,1024,2048\}$. }}
\label{histogram}
\end{center}
\end{figure}

In the lower part of Table \ref{vglmitAUE1} we summarize the simulated power of the test \eqref{test} and the method proposed by \cite{aue2009} for various configurations of the model \eqref{model4}. It can be observed that the new method significantly outperforms the method of \cite{aue2009} for models with more than one break for small sample sizes. For the configurations featuring one break point we observe that the new method performs similarly for large sample sizes, while providing slightly better performance for small sample sizes. These observations are remarkable since model \eqref{model4} defines a process with switching variance and the test of \cite{aue2009} is particularly designed for detecting this structure. Since our approach is able to detect a much wider class of alternatives one might expect a loss in power compared to a procedure which is specifically constructed to test for such changes. Surprisingly, this does not seem to be the case. \\
Next we illustrate the performance of the procedure for the localization of the break points. For this purpose we simulated data from model \eqref{beispielbild} and applied the new method to detect and estimate the location of the change points. We repeated this 100 times and for each trial we obtained an estimate $\hat{\mathfrak{b}}=(\hat b_1,...,\hat b_{\hat K})$ of the location $\mathfrak{b}=(b_1,...,b_K)$. The histograms in Figure \ref{histogram} present the obtained empirical distribution of the estimated break points for different sample sizes $T$. We observe that all histograms are centered at the true break points and that the variation decreases with increasing sample size.

\begin{sidewaysfigure}
\includegraphics[width=0.95\textheight, height=120mm]{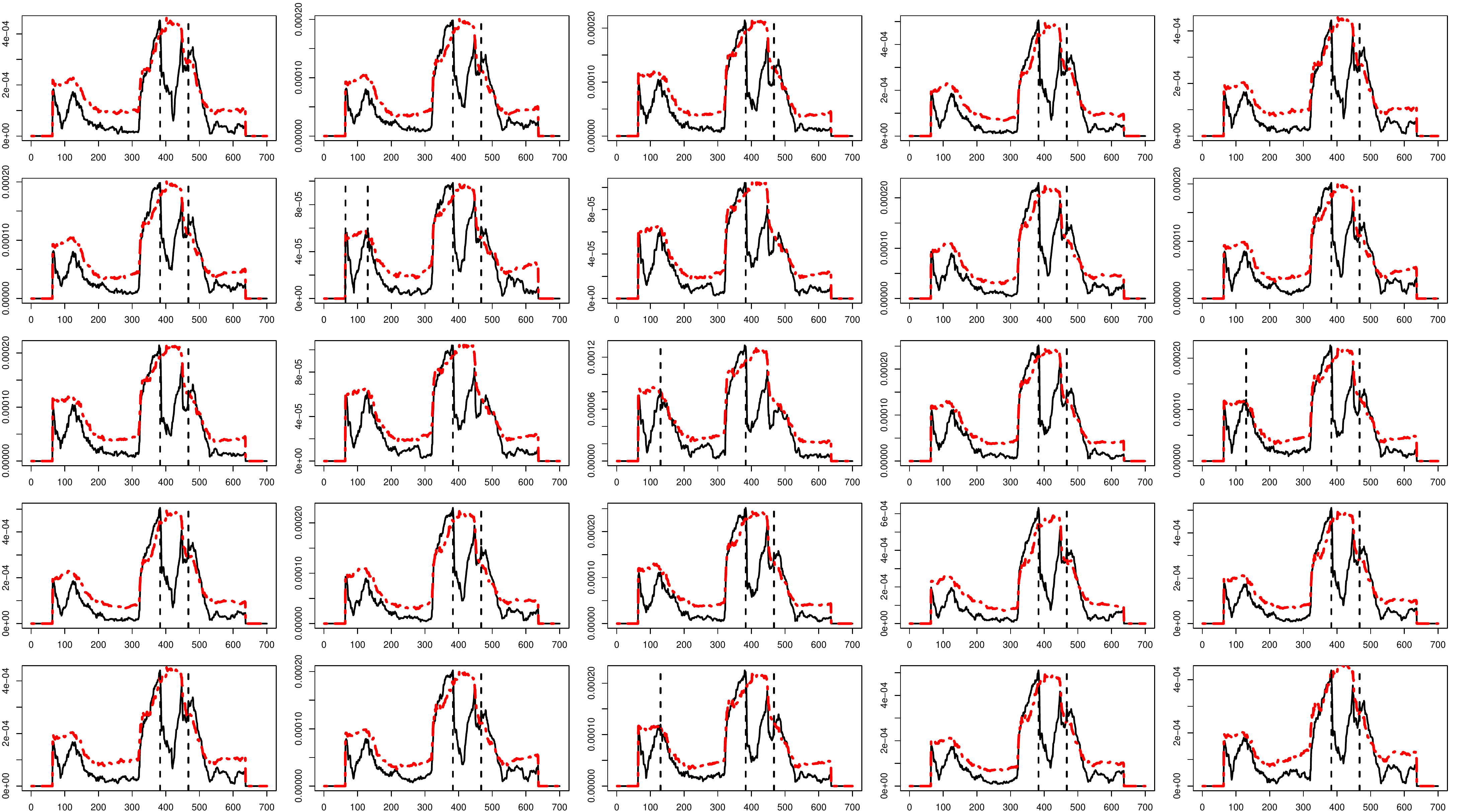}
\caption{\textit{The functions $v \mapsto N^\gamma \sup\limits_{\omega \in [0,1]}\big|[\hat{\mathbf{D}}_T(v,\omega)]_{a,b}\big|$ (solid lines) and $v \mapsto \boldsymbol{\epsilon}_{T,a,b}(v)$ (broken lines) for $a,b=1,...,5$ for the daily log returns (November 2010 - November 2012) corresponding to the five sector ETFs. In each plot the dashed vertical lines mark the estimated break points.}}
\label{poud2}
\end{sidewaysfigure}

\begin{figure}
\begin{center}
\includegraphics[width=160mm,height=150mm]{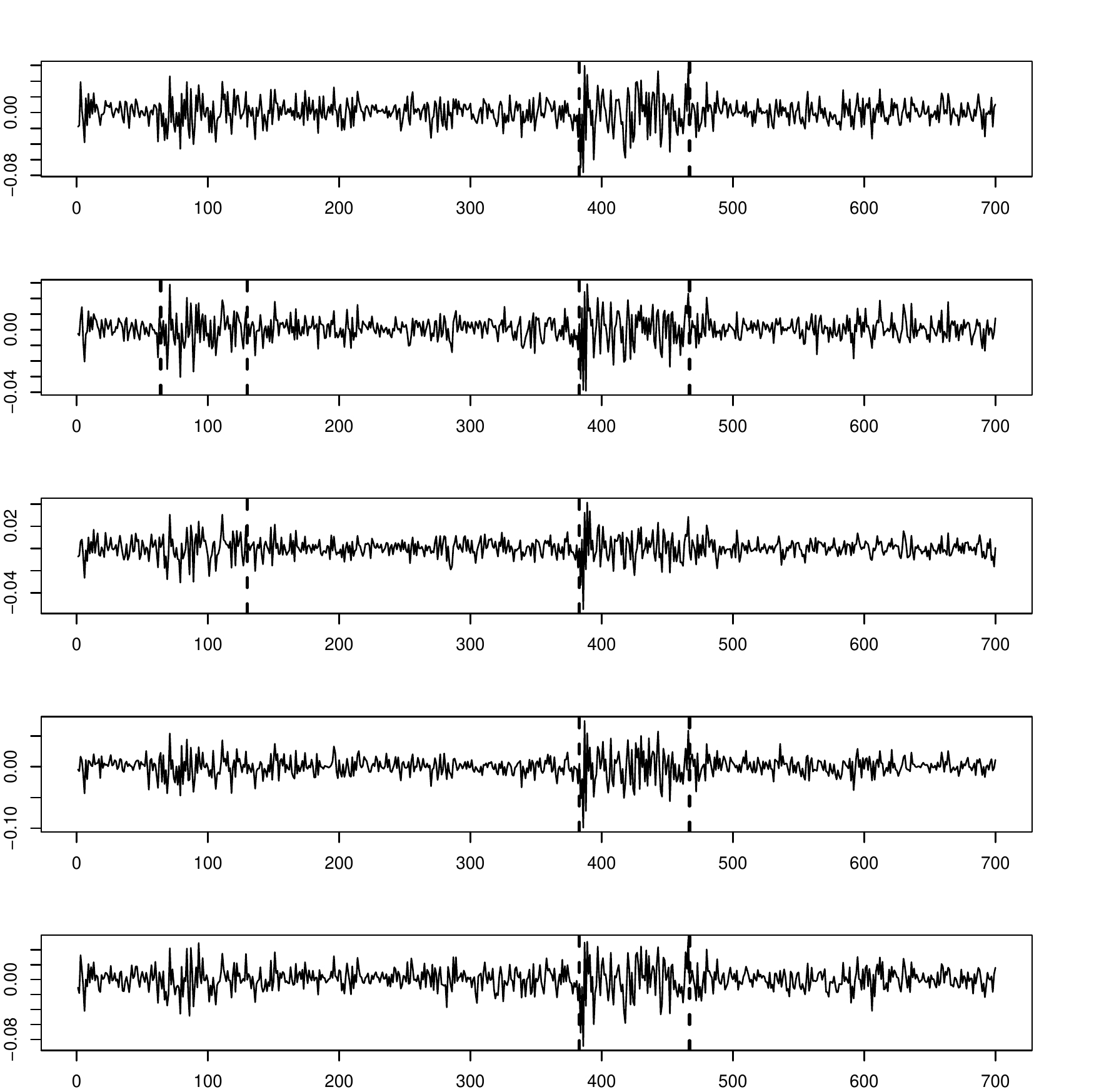}
\caption{\textit{The log returns of the sector ETFs. In each plot the vertical lines mark the estimated break points of the {univariate} process. }}
\label{poud3}
\end{center}
\end{figure}

We conclude this section with an illustration of the procedures in a multivariate data example. For this purpose we investigate five stock ETFs (Exchange Traded Funds), namely Materials Select Sector SPDR (XLB), Consumer Staples Select Sector SPDR (XLP), Utilities Select Sector SPDR (XLU), Financial Select Sector SPDR (XLF) and Energy Select Sector SPDR (XLE). In Figure \ref{poud3} the five dimensional time series $\{\boldsymbol{X_{t}}\}_{t=1,...,700}$ of the daily log returns corresponding to the five ETFs between February 2010 and November 2012 is displayed. The new method yields the following results: The test \eqref{test} for structural breaks is highly significant and rejects the null hypothesis with a $p$-value of zero. In a second step we use the detection procedure with $N=64$ and detect four break points at dates 29/04/2010, 03/08/2010, 03/08/2011, and 01/12/2011. We then also applied the refined analysis in order to identify the components where the change point is present. The results are displayed in Figure \ref{poud2} and \ref{poud3}. 
Note  that at the third break point 03/08/2011 each component $N^{\gamma}\sup\limits_{\omega\in [0,1]}\big|[\hat{\mathbf{D}}_T(\hat b_i,\omega)]_{a,b}\big|$ surpasses the respective critical threshold $\boldsymbol{\epsilon}_{a,b,T}(\hat b_i)$ and we conclude that each component of the spectral density exhibits a change point. On the other hand, for the other three break points only some entries of this matrix exceed the threshold sequence and we therefore detect a structural break only in some components of the spectral density matrix.

{\bf Acknowledgements}
This work has been supported in part by the
Collaborative Research Center ``Statistical modeling of nonlinear
dynamic processes'' (SFB 823, Teilprojekt  C1) of the German Research Foundation
(DFG).

\vspace{-.5cm}
\bibliographystyle{apalike}

\section{Appendix I: Sketch of the proofs}  \label{sec7}
\def\theequation{7.\arabic{equation}}
\setcounter{equation}{0}

 For the sake of brevity and in order to improve the readability we 
restrict ourselves to the main ideas of the proofs in this section and present all technical details in an additional Appendix thereafter. \\

{\bf Proof of Theorem \ref{hauptsatz}}
For notational convenience we restrict ourselves to the case $d=1$, since the more general case is treated completely analogously using linearity arguments and the independence of the components of $\boldsymbol{Z}_t$. Throughout this chapter $C$ denotes a universal constant, which does not depend on the sample size and can vary from line to line in the calculations.\\
{\it Proof of part a):} For the proof of
\eqref{austheorem} it is sufficient to show the following two claims: [see Theorem 1.5.4 and 1.5.7 in \cite{wellnervandervaart}]:
\begin{itemize}
\item[(1)] For every $k\in\N$ and $y_1:=(v_1,\omega_1),...,y_k:=(v_k,\omega_k) \in [0,1]$ we have
\vspace{-.5cm}
\begin{align}
\label{finitedimensional}
\sqrt{N}[\hat D_T(y_1),...,\hat D_T(y_k)]\Rightarrow[G(y_1),...,G(y_k)].
\end{align}
\item[(2)] For every $\eta, \epsilon>0$ there exists a $\delta>0$ such that
\vspace{-.5cm}
\begin{align}
\label{gleichstetigkeit}
\lim_{T\rightarrow \infty} P\Big( \sup_{(y_1, y_2) \in [0,1]^2: d_2(y_1,y_2) <\delta } N^{1/2} |\hat D_T(y_1)-\hat D_T(y_2)| > \eta \Big)< \epsilon ,
\end{align} 
where $d_2(y_1,y_2)$ denotes the euclidean distance between $y_1=(v_1,\omega_1)$ and $y_2=(v_2,\omega_2)$.
\end{itemize}

{\it Proof of \eqref{finitedimensional}:} The assertion follows if we are able to show that, for each $k\in\N$ and each $y_1=(v_1,\omega_1),...,y_k=(v_k,\omega_k)$, all cumulants of the random vector $[\sqrt{N}\hat D_T(y_i)]_{i=1,...,k}$ converge to the corresponding cumulants of the vector $[G(y_i)]_{i=1,...,k}$. We thus have to show for all $y=(v,\omega),y_1,...,y_l\in[0,1]^2$:
\begin{itemize}
\item[(i)]$\E(\sqrt{N}\hat D_T(y))=o(1)$.
\item[(ii)]$\text{Cov}(\sqrt{N}\hat D_T(y_1),\sqrt{N}\hat D_T(y_2))=\text{Cov}(G(y_1),G(y_2))+o(1)$. 
\item[(iii)]$\text{cum}(\sqrt{N}\hat D_T(y_1),...,\sqrt{N}\hat D_T(y_l))=o(1)$.
\end{itemize}
For this purpose we define for $y=(v,\omega)$
\begin{align}\label{PHI}
\phi_{y,T}(j,\lambda):=&1_{[0,\frac{2\pi\lfloor \omega N/2 \rfloor}{N}]}(\lambda)\Big[1_{\{\lfloor u(v,T) T\rfloor+   N/2   \}}(j)-1_{\{\lfloor u(v,T) T\rfloor-   N/2 \}}(j)\Big],
\end{align}
where the function $u$ is defined by
$ u(v,T) =v $ if  $N/T\leq v \leq 1-N/T$, $  N/T $ {if } $ v < N/T$,  and 
$1-N/T$   {if } $v>1-N/T$. 
This notation implies the representation
\vspace{-.5cm}
\begin{align}\label{REPRESENTATION}
\hat D_T(y)&=\frac{1}{N}\sum_{j=1}^{T}\sum_{k=1}^{   N/2   }\phi_{y,T}(j,\lambda_k)I_N(\frac{j}{T},\lambda_k).
\end{align}
 Additionally, we set
$
D_{N,T}(y):=\frac{T}{N} (\int_{0}^{\omega\pi}\int_{v}^{v+N/T}f(u,\lambda)dud\lambda-\int_{0}^{\omega\pi}\int_{v-N/T}^vf(u,\lambda)dud\lambda),
$
and start with a proof of
\vspace{-.5cm}
\be \label{konsistenzew}
\E\Big(\sqrt{N}\big(\hat D_T(y)-D_{N,T}(y)\big)\Big)=o(1),
\ee
from which (i)  follows directly because of $D_{N,T}(y) \equiv 0$ under the null hypothesis [we treat the more general case since \eqref{konsistenzew} is also required in the proof of part b)]. By writing $\psi_l(t/T):=\boldsymbol{\Psi}_l(t/T)$ we get $(y=(v,\omega))$
\be
&&\E\big(\hat D_T(y)\big)=\frac{1}{2\pi N^2}\sum_{j=1}^T\sum_{k=1}^{N/2}\phi_{y,T}(j,\lambda_k)\sum_{p,q=0}^{N-1}e^{-i\lambda_k(p-q)}\E(X_{j-\frac{N}{2}+1+p}X_{j-\frac{N}{2}+1+q})\nonumber\\
&=&\frac{1}{2\pi N^2}\sum_{j=1}^T\sum_{k=1}^{N/2}\phi_{y,T}(j,\lambda_k)\sum_{l,m=0}^{\infty}\psi_l(\frac{j-\frac{N}{2}+1+p}{T})\psi_m(\frac{j-\frac{N}{2}+1+q}{T})\sum_{p,q=0}^{N-1}e^{-i\lambda_k(p-q)}\nonumber\\
&&\times\E(Z_{j-\frac{N}{2}+1+p-l}Z_{j-\frac{N}{2}+1+q-m})\label{gl1}
\ee
and by using the identity $\E(Z_iZ_j)=\delta_{ij}$ [here and throughout this paper $\delta_{ij}$ denotes the Kronecker symbol] we obtain that the restriction $q=p-l+m$ has to hold such that the respective summands in \eqref{gl1} do not vanish. If we furthermore define $A_{T,1}(v):=\{\lfloor vT\rfloor-N/2,\lfloor vT\rfloor+N/2\}$, we obtain for \eqref{gl1}
\be
&&\frac{1}{2\pi N^2}\sum_{j\in A_{T,1}(v)}\sum_{k=1}^{N/2}\phi_{y,T}(j,\lambda_k)\sum_{l,m=0}^{\infty}\sum_{\substack{p=0\\0\leq p-l+m\leq N-1}}^{N-1}\nonumber\\
&&\times\psi_l(\frac{j-\frac{N}{2}+1+p}{T})\psi_m(\frac{j-\frac{N}{2}+1+p-l+m}{T}) e^{-i\lambda_k(m-l)}\nonumber\\
&=&\frac{1}{2\pi N^2}\sum_{j\in A_{T,1}(v)}\sum_{k=1}^{N/2}\phi_{y,T}(j,\lambda_k)\sum_{l,m=0}^{\infty}\sum_{\substack{p=0\\0\leq p-l+m\leq N-1}}^{N-1} \Big \{ \psi_l(\frac{j-\frac{N}{2}+1+p}{T})\psi_m(\frac{j-\frac{N}{2}+1+p}{T})\nonumber\\
&+&\psi_l(\frac{j-\frac{N}{2}+1+p}{T})\Big(\psi_m(\frac{j-\frac{N}{2}+1+p-l+m}{T})-\psi_m(\frac{j-\frac{N}{2}+1+p}{T})\Big) \Big \} e^{-i\lambda_k(m-l)}\nonumber\\
&=:&E_{1,T}(y)+E_{2,T}(y)\label{gl2},
\ee
where $E_{1,T}$ and $E_{2,T}$ are defined in an obvious manner.
The assertion now follows from
\begin{itemize}
\item[(ia)]$E_{1,T}(y)=D_{N,T}(y)+O(1/N)$.
\item[(ib)]$E_{2,T}(y)=o(1/\sqrt{N})$
\end{itemize}
which are established in Appendix II. The proof of (ii) and (iii) are also given in Appendix II.  \hfill $\Box$

{\it Proof of  \eqref{gleichstetigkeit}:} Observing
\vspace{-.5cm}
\begin{align}
&\sqrt{N}\hat D_T(y)=\frac{1}{\sqrt{N}}\sum_{k=1}^{\lfloor\omega N/2\rfloor}I_N(u(v,T)+N/(2T),\lambda_k)-\frac{1}{\sqrt{N}}\sum_{k=1}^{\lfloor\omega N/2\rfloor}I_N(u(v,T)-N/(2T),\lambda_k)\nonumber\\
&=:\sqrt{N}\hat D^{(1)}_T(y)-\sqrt{N}\hat D^{(2)}_T(y),\label{decomposition1}
\end{align}
it is sufficient to show asymptotic stochastic equicontinuity for the processes \linebreak$\{\sqrt{N}\hat D^{(1)}_T(y)\}_{y\in[0,1]^2}$ and $\{\sqrt{N}\hat D^{(2)}_T(y)\}_{y\in[0,1]^2}$. We only present the proof for the first summand in \eqref{decomposition1} and note that stochastic equicontinuity for the second term can be shown analogously. With the notation
$
\phi^{(1)}_{y,T}(j,\lambda):=1_{[0,\frac{2\pi\lfloor \omega N/2 \rfloor}{N}]}(\lambda)1_{\{\lfloor u(v,T)\times T\rfloor+   N/2   \}}(j),
$
we obtain the representation
$
\hat D^{(1)}_T(y)=\frac{1}{N}\sum_{j=1}^{T}\sum_{k=1}^{   N/2   }\phi^{(1)}_{y,T}(j,\lambda_k)I_N(j/T,\lambda_k).
$
For $y_i=(v_i,\omega_i)$ ($i=1,2$) we define a semi-metric
$
d_T(y_1,y_2):=\sqrt{|\omega_2-\omega_1|+|\lfloor v_1T \rfloor - \lfloor v_2 T \rfloor |/N}
$
on the set $\mathcal{P}_T:=\{0,1/T,2/T,...,1-N/T\}\times\{1/N,...,1-1/N,1\}$. This yields
\bea
\Delta_{\delta,\eta}:=\mathbb{P}(\sup\limits_{\substack{y_i\in[0,1]^2\\ d_2(y_1,y_2)<\delta}}\sqrt{N}|&\hat D^{(1)}_T(y_1)-\hat D^{(1)}_T(y_2) |>\eta)=\mathbb{P}(\sup\limits_{\substack{y_i\in \mathcal{P}_T\\ d_2(y_1,y_2)<\delta}}\sqrt{N}|\hat D^{(1)}_T(y_1)-\hat D^{(1)}_T(y_2) |>\eta),
\eea

and it is easy to verify that, for a fixed $\delta>0$, there exists a $\delta'>0$ [with $\delta'(\delta)\rightarrow 0$ as $\delta\rightarrow 0$] such that
\vspace{-.5cm}
\begin{align}
\Delta_{\delta,\eta}\leq \mathbb{P}(\sup\limits_{\substack{y_i\in \mathcal{P}_T:d_T(y_1,y_2)<\delta'}}\sqrt{N}|\hat D^{(1)}_T(y_1)-\hat D^{(1)}_T(y_2) |>\eta)\label{ungl3}
\end{align}
is fulfilled if $T$ is sufficiently large. So it suffices to prove that the probability on the right hand side of \eqref{ungl3} can be made arbitrarily small if $T$ is sufficiently large. For this purpose let $C(u,d_T,\mathcal{P}_T)$ denote the covering number of $\mathcal{P}_T$ with respect to the semi-metric $d_T(\cdot,\cdot)$ and define the corresponding covering integral of $\mathcal{P}_T$ by
$
J_T(\kappa):=\int_0^{\kappa}[\log\Big(\frac{48C(u,d_T,\mathcal{P}_T)^2}{u}\Big)]^2 du.
$
In Appendix II we establish the assertions 
\vspace {-.5cm}
\begin{eqnarray}
&& \lim\limits_{\kappa \rightarrow 0}\lim\limits_{T \rightarrow \infty}J_T(\kappa)=0, \label{cond0}
 \\
 && 
\label{cond1}
\mathbb{E}\Big(N^{k/2}\big(\hat{D}^{(1)}_T(y_1)-\hat{D}^{(1)}_T(y_2) \big)^k\Big)\leq (2k)!C^kd_T\big(y_1,y_2\big)^k
\end{eqnarray}
for  a constant $C \in \R^+$  and  all $y_1,y_2\in [0,1]^2$ and even integers $k\in\N$.
Therefore it follows, by similar arguments as given in \cite{dahlhaus1988}, that
\vspace {-.5cm}
\bea
\mathbb{P}(\sup\limits_{\substack{y_i\in \mathcal{P}_T\\d_T(y_1,y_2)<\delta'}}\sqrt{N}|\hat D^{(1)}_T(y_1)-\hat D^{(1)}_T(y_2) |>\eta) < \epsilon
\eea
for $T$ sufficiently large and sufficiently small $\delta'=\delta'(\delta)>0$, which proves stochastic equicontinuity.  \hfill $\Box$

{\it Proof of part b):} Under the alternative,  there exist  for all $r \in \{1,...,K\}$ an $\omega_r$ such that $|D(b_r,\omega_r)|>0$. Note that, the proof of part a) does not rely on the property that the functions $\psi_l(u)$ are constant. In fact, only the treatment of the expectation [i.e. part (1) (i)] is slightly easier if $\psi_l(u)=\psi_l$. However, in this case we proved the more general claim \eqref{konsistenzew}. So, by following the proofs of $(ii)$ and $(iii)$ in the proof of part a) and employing \eqref{konsistenzew}, we obtain 
$
N^{1/2}\Big|\Big|\hat D_{T}(y) -  D_{N,T}(y)\Big|\Big|_{\infty}=O_P(1),
$
 which directly yields the assertion. $\hfill \Box$
 \bigskip

{\bf Proof of Theorem \ref{hauptsatz2}:} See Appendix II.

{\bf Proof of Theorem \ref{theorembootstrap}}
As in the previous proof we restrict ourselves without loss of generality to
the case $d = 1$. Furthermore we suppress the argument $T$, when referring to the sequence $p = p(T)$. Because of Assumption \ref{annahmenbootstrap}, we obtain the $MA(\infty)$ representation
$
X_{t,T}^*=\sum_{l=0}^\infty \hat\psi_l^{AR}(p)Z_{t-l}^*
$
for $T$ and $p(T)$ sufficiently large [see Section $3$ of \cite{detprevet2011b} for more details]. Note that this representation corresponds to a process without structural breaks and that in the proof of Theorem \ref{hauptsatz}a)  all error terms can be bounded by
\begin{align*}
\frac{(\sum_{m=0}^\infty |\psi_m|)^{q_1}(\sum_{l=0}^\infty l |\psi_l|)^{q_2}}{N}=O(1/N)
\end{align*}
where $q_1,q_2 \in \N$ and  the equality is a consequence of \eqref{summ1}. So the proof of Theorem \ref{theorembootstrap} follows in the same way as the proof of Theorem \ref{hauptsatz} a) if we show that the (now random) errors terms are of order $O_P(1/N)$. However, it was shown in Theorem 3.2 of \cite{detprevet2011b} that 
$
(\sum_{m=0}^\infty |\hat \psi_m^{AR}(p)|)^{q_1}( \sum_{l=0}^\infty l |\hat \psi_l^{AR}(p)|)^{q_2}=O_P(1)
$
holds for all $q_1,q_2 \in \N$, which directly yields the claim. $\hfill \Box$
\medskip

{\bf Proof of Theorem \ref{theorembootstrap2}:}  see Appendix II.

{\bf Proof of Theorem \ref{detection}}
For a Proof of part a) note that  {$\sup\limits_{\omega\in[0,1]}|[\boldsymbol{D}_{N,T}(v,\omega)]_{a,b}|=0$} for $v\in\bar{\mathcal{I}}_{T,a,b}(b_1,...,b_K)$, and observe \eqref{propertydetection}.
This  yields for all $a,b \in \{1,...,d\}$ and sufficiently large $N,T$
\vspace{-.35cm}
\begin{align*}
&\mathbb{P}\Big(\bigcup\limits_{v \in \bar{\mathcal{I}}_{T,a,b}(b_1,...,b_K)}\Big\{N^{\gamma}\sup\limits_{\omega \in [0,1]}|[\boldsymbol{\hat{D}}_T(v,\omega)]_{a,b}|>\boldsymbol{\epsilon}_{T,a,b}(v)\Big\}\Big)\\
\leq
&\mathbb{P}\Big(N^{\gamma}\sup\limits_{v \in \bar{\mathcal{I}}_{T,a,b}(b_1,...,b_K)}\sup\limits_{\omega \in [0,1]}|[\boldsymbol{\hat{D}}_T(v,\omega)]_{a,b}-[\boldsymbol{D}_{N,T}(v,\omega)]_{a,b}|>C/2\Big)\\
\leq&\mathbb{P}\Big(N^{\gamma}\sup\limits_{v \in [0,1]}\sup\limits_{\omega \in [0,1]}||\boldsymbol{\hat{D}}_T(v,\omega)-\boldsymbol{D}_{N,T}(v,\omega)||_{\infty}>C/2\Big)\xrightarrow{\text{ }T \rightarrow \infty\text{ }} 0.
\end{align*}
where we used  Theorem \ref{hauptsatz} (if assumption (2) is fulfilled)  and \ref{hauptsatz2} (if assumption (1) is fulfilled).
{Part b)} is a direct consequence of Theorem \ref{hauptsatz}b) and \ref{hauptsatz2}b), which imply 
for  $r=1,...,K$ and  $(a,b) \in B(b_r)$ : 
$\mathbb{P}(\sup\limits_{\omega\in[0,1]}N^{\gamma}|[\boldsymbol{\hat{D}}_T(b_r,\omega)]_{a,b}|>\boldsymbol{\epsilon}_{Ta,b}(b_r) )\rightarrow 1
$.  \hfill $\Box$

\section{Appendix II: technical details}  \label{sec8}
\def\theequation{8.\arabic{equation}}
\setcounter{equation}{0}

\subsection{Proof of (ia) and (ib) in the proof of part (a) Theorem \ref{hauptsatz}}
For a proof of (ia), we note that by \eqref{summ1} we can drop the restriction in the summation with respect to $p$ by making an error of order $O(1/N)$. Thus we obtain for $E_{1,T}$
\be
&&\frac{1}{2\pi N^2}\sum_{j\in A_{T,1}(v)}\sum_{k=1}^{N/2}\phi_{y,T}(j,\lambda_k)\sum_{l,m=0}^{\infty}\sum_{p=1}^{N-1}\psi_l(\frac{j-\frac{N}{2}+1+p}{T})\psi_m(\frac{j-\frac{N}{2}+1+p}{T}) e^{-i\lambda_k(m-l)}\nonumber\\
&&+O(1/N)\nonumber\\
&=&\frac{1}{N}\sum_{j\in A_{T,1}(v)}\sum_{k=1}^{N/2}\phi_{y,T}(j,\lambda_k)\frac{1}{2\pi }\sum_{l,m=0}^{\infty}\frac{T}{N}\int_{j/T-\frac{N}{2T}+\frac{1}{T}}^{j/T+\frac{N}{2T}}\psi_l(u)\psi_m(u)du\quad e^{-i\lambda_k(m-l)}du+O(1/N)\nonumber \\
&=&D_{N,T}(y)+O(1/N),
\ee
where the second equality follows with the piecewise constancy of the functions $\psi_l(u)$. For (ib)  we get
\begin{align*}
&|E_{2,T}(y)|\leq C\sum_{j\in A_{T,1}(v)}\sum_{\substack{l,m=0\\|l-m|\leq N }}^{\infty}\Big|\int_{j/T-\frac{N}{2T}+\frac{1}{T}}^{j/T+\frac{N}{2T}}\psi_l(u)\big(\psi_m(u+\frac{m-l}{T})-\psi_m(u)\big)du\Big|+O(1/N),
\end{align*}

 where the summation can be restricted to indices satisfying $|l-m|\leq N$, which follows directly from the restriction on $p$. By employing that
\begin{equation*}
\int_a^bf(x)g(x+y)dx-\int_a^bf(x)g(x)dx\leq C |y|\sup\limits_{z\in[a-|y|,b+|y|]}|f(z)|\sup\limits_{z_1,z_2\in[a-|y|,b+|y|]}|g(z_1)-g(z_2)|,
\end{equation*}

holds for all  piecewise constant  functions $f$, $g$ exhibiting only finitely many points of discontinuity on the interval $[a-|y|,b+|y|]$, we obtain
\begin{align*}
|E_{2,T}(y)\leq|&C\sum_{j\in A_{T,1}(v)}\sum_{\substack{l,m=0\\|l-m| \leq N}}^{\infty}\frac{|m-l|}{T}\sup\limits_{u\in[0,1]}|\psi_l(u)|\sup\limits_{u\in[0,1]}|\psi_m(u)|+O(1/N)=o(1/\sqrt{N}).
\end{align*}

\subsection{Proof of (ii) and (iii) in the proof of part (a) of Theorem \ref{hauptsatz}}

For the proof of part (ii) note that, under the null hypothesis of no structural breaks, the quantities $\psi_l(t/T)$ do not depend on the rescaled time  $t/T$ and we can work with the representation
\bea
X_{t,T}=\sum_{l=0}^{\infty}\psi_lZ_{t-l}.
\eea
Without loss of generality we assume that $v_1\leq v_2$ holds and define the set
$$A_{T,2}(v_1,v_2):=\{\lfloor v_1T\rfloor-   N/2   ,\lfloor v_1T\rfloor+   N/2\}\times   \{\lfloor v_2T\rfloor-   N/2   ,\lfloor v_2T\rfloor+   N/2   \}.$$
We then obtain for $\omega_1,\omega_2\in[0,1]$
\begin{align}
&\Cov(\sqrt{N}\hat D_T(y_1),\sqrt{N}\hat D_T(y_2))\nonumber\\
=&\frac{1}{N}\sum_{(j_1,j_2)\in A_{T,2}(v_1,v_2)}\sum_{k_1,k_2=1}^{N/2}\phi_{y_1,T}(j_1,\lambda_{k_1})\phi_{y_2,T}(j_2,\lambda_{k_2})\cum(I_N(\frac{j_1}{T},\lambda_{k_1}),I_N(\frac{j_2}{T},\lambda_{k_2}))\nonumber\\
=&\frac{1}{(2\pi)^2 N^3}\sum_{(j_1,j_2)\in A_{T,2}(v_1,v_2)}\sum_{k_1,k_2=1}^{N/2}\phi_{y_1,T}(j_1,\lambda_{k_1})\phi_{y_2,T}(j_2,\lambda_{k_2})\sum_{p_1,p_2=0}^{N-1}\sum_{q_1,q_2=0}^{N-1}\sum_{l,m,n,o=0}^{\infty}\psi_l\psi_m\nonumber\\
&\times\psi_n\psi_o e^{-i\lambda_{k_1}(p_1-q_1)} e^{-i\lambda_{k_2}(p_2-q_2)}\cum(Z_{j_1+p_1+1-l}Z_{j_1+q_1+1-m},Z_{j_2+p_2+1-n}Z_{j_2+q_2+1-o}) \nonumber\\
=&:\sum_{(j_1,j_2)\in A_{T,2}(v_1,v_2)}B_T(j_1,j_2)\label{cov0},
\end{align}
where $B_T(j_1,j_2)$ is defined in an obvious manner. Using
\bea
\cum(Z_{a}Z_b,Z_cZ_d)=\cum(Z_a,Z_d)\cum(Z_b,Z_c)+\cum(Z_a,Z_c)\cum(Z_b,Z_d)
\eea
[see Theorem 2.3.2 in \cite{brillinger1981}] we can split $B_T(j_1,j_2)$  into two parts, and the independence of the innovations $Z_t$ yields for the first term
\begin{align}
&\frac{1}{(2\pi)^2 N^3}\sum_{k_1,k_2=1}^{N/2}\phi_{y_1,T}(j_1,\lambda_{k_1})\phi_{y_2,T}(j_2,\lambda_{k_2})\sum_{l,m,n,o=0}^{\infty}\psi_l\psi_m\psi_n\psi_o\nonumber\\
&\times \sum_{\substack{p_1,p_2=0\\0\leq p_1-l+o+j_1-j_2\leq N-1\\ 0\leq p_2-n+m+j_2-j_1\leq N-1}}^{N-1} e^{-i\lambda_{k_1}(p_1-p_2+n-m-j_2+j_1)}e^{-i\lambda_{k_2}(p_2-p_1+l-o-j_1+j_2)}\nonumber\\
=&\frac{1}{(2\pi)^2 N^3}\sum_{k_1,k_2=1}^{N/2}\phi_{y_1,T}(j_1,\lambda_{k_1})\phi_{y_2,T}(j_2,\lambda_{k_2})\sum_{l,m,n,o=0}^{\infty}\psi_l\psi_m\psi_n\psi_o e^{-i\lambda_{k_1}(n-m)} e^{-i\lambda_{k_2}(l-o)}\nonumber\\
&\times\sum_{\substack{p_1,p_2=0\\0\leq p_1-l+o+j_1-j_2\leq N-1\\ 0\leq p_2-n+m+j_2-j_1\leq N-1}}^{N-1} e^{-i(\lambda_{k_1}-\lambda_{k_2})(p_1-p_2-j_2+j_1)}=:V_{k_1=k_2}+V_{k_1\neq k_2}, \nonumber
\end{align}
where $V_{k_1=k_2}$ and $V_{k_1\neq k_2}$ denote the summation of all terms with $k_1=k_2$ and $k_1\neq k_2$ respectively [note that the restrictions $p_1-l+o+j_1-j_2=q_2$ and $p_2-n+m+j_2-j_1=q_1$ follow by the independence of the innovations $Z_t$]. For the first term we obtain
\begin{align}
V_{k_1=k_2}=&\frac{1}{(2\pi)^2 N^3}\sum_{k=1}^{N/2}\phi_{y_1,T}(j_1,\lambda_{k})\phi_{y_2,T}(j_2,\lambda_{k})\sum_{l,m,n,o=0}^{\infty}\psi_l\psi_m\psi_n\psi_o e^{-i\lambda_{k}(n-m)}e^{-i\lambda_{k}(l-o)}\nonumber\\
&\times\max(N-1-|j_1-j_2-l+o|,0)\max(N-1-|j_2-j_1+m-n|,0).\label{cov3}
\end{align}
By applying the summability condition \eqref{summ1} we obtain that $V_{k_1=k_2}$ is of order $O(\frac{1}{N})$ if $|j_1-j_2|>N$, and that for $\Delta:=|j_1-j_2|<N$,  \eqref{cov3} is the same as
\begin{align*}
\frac{1}{N}\big(1-\frac{\Delta}{N}\big)^2\sum_{k=1}^{N/2}\phi_{y_1,T}(j_1,\lambda_{k})\phi_{y_2,T}(j_2,\lambda_{k})f^2(\lambda_k)+O(\frac{1}{N}).
\end{align*}
For the quantity $V_{k_1\neq k_2}$ we get by simple calculations and an application of \eqref{summ1}
\begin{align}
V_{k_1\neq k_2}=&\frac{1}{(2\pi)^2N^3}\sum_{\substack{k_1,k_2=1\\k_1\neq k_2}}^{N/2}\phi_{y_1,T}(j_1,\lambda_{k_1})\phi_{y_2,T}(j_2,\lambda_{k_2})\sum_{l,m,n,o=0}^{\infty}\psi_l\psi_m\psi_n\psi_oe^{-i\lambda_{k_1}(n-m)}e^{-i\lambda_{k_2}(l-o)} \nonumber\\
&\times\sum_{\substack{p_1,p_2=0\\0\leq p_1+j_1-j_2\leq N-1\\ 0\leq p_2+j_2-j_1\leq N-1}}^{N-1}e^{-i(\lambda_{k_1}-\lambda_{k_2})(p_1-p_2-j_2+j_1)}+O(\frac{1}{N})\nonumber
\end{align}
[i.e. we can omit the $l,m,n,o$ in the restrictions on $p_i$ at the cost of a term of order $O(\frac{1}{N})$]. In the following let $a_T$ denote a sequence satisfying $a_T\rightarrow\infty$, $a_T/N\rightarrow 0$ and $N^2 /a_T^3 \rightarrow 0$. It is straightforward to verify
\bea
\sum_{\substack{p_1,p_2=0\\0\leq p_1+j_1-j_2\leq N-1\\ 0\leq p_2+j_2-j_1\leq N-1}}^{N-1}e^{-i(\lambda_{k_1}-\lambda_{k_2})(p_1-p_2-j_2+j_1)}=\Big|\sum_{p=0}^{N-1-\Delta}e^{-i(\lambda_{k_1}-\lambda_{k_2})p}\Big|^2=\Big|\frac{1-e^{i\frac{2\pi(k_1-k_2)}{N}\Delta}}{1-e^{-i\frac{2\pi (k_1-k_2)}{N}}}\Big|^2=\Big|\frac{\sin(\frac{\pi (k_1-k_2)}{N}\Delta)}{\sin(\frac{\pi (k_1-k_2)}{N})}\Big|^2,
\eea
which yields
\bea
&&V_{k_1\neq k_2}=\frac{1}{N^3}\sum_{\substack{k_1,k_2=1\\k_1\neq k_2}}^{N/2}\phi_{y_1,T}(j_1,\lambda_{k_1})\phi_{y_2,T}(j_2,\lambda_{k_2})f(\lambda_{k_1})f(\lambda_{k_2})\Big|\frac{\sin(\frac{\pi (k_1-k_2)}{N}\Delta)}{\sin(\frac{\pi (k_1-k_2)}{N})}\Big|^2+O(\frac{1}{N}) \\
&&=\frac{1}{N^3}\sum_{\substack{k_1,k_2=1\\k_1\neq k_2\\|k_1-k_2|\leq a_T}}^{N/2}\phi_{y_1,T}(j_1,\lambda_{k_1})\phi_{y_2,T}(j_2,\lambda_{k_2})f(\lambda_{k_1})f(\lambda_{k_2})\Big|\frac{\sin(\frac{\pi (k_1-k_2)}{N}\Delta)}{\sin(\frac{\pi (k_1-k_2)}{N})}\Big|^2+O(\frac{1}{N}+\frac{1}{a_T^{1-\tau}})
\eea
for any $\tau>0$, because of
\begin{align*}
\frac{1}{N}\sum_{\substack{k_1,k_2=1\\|k_1-k_2|>a_T}}^{N/2}\frac{1}{N^2}\frac{1}{\sin^2(\frac{\pi(k_1-k_2)}{N})}\leq\frac{C}{N}\sum_{\substack{k_1,k_2=1\\|k_1-k_2|>a_T}}^{N/2}\frac{1}{(k_1-k_2)^2}\leq\frac{1}{a_T^{1-\tau}}.
\end{align*}
By  applying $\sin(x)=x+O(\epsilon^3)$ for $x\in[0,\epsilon]$, a Taylor expansion and $\sum_{k=1}^n 1/k \sim \log(n)$ we obtain, that $V_{k_1\neq k_2}$ is the same as
\begin{align}
&\frac{1}{N}\sum_{\substack{k_1,k_2=1\\k_1\neq k_2\\|k_1-k_2|\leq a_T}}^{N/2}\phi_{y_1,T}(j_1,\lambda_{k_1})\phi_{y_2,T}(j_2,\lambda_{k_1})f^2(\lambda_{k_1})\Big|\frac{\sin(\frac{\pi (k_1-k_2)}{N}\Delta)}{\pi (k_1-k_2)}\Big|^2+O(\frac{1}{N}+\frac{1}{a_T^{1-\tau}}+\frac{N^2}{a_T^3}). \nonumber
\end{align}
This expression  is equal to
\begin{align}
&\frac{1}{N}\sum_{k_1=a_T}^{N/2-a_T}\phi_{y_1,T}(j_1,\lambda_{k_1})\phi_{y_2,T}(j_2,\lambda_{k_1})f^2(\lambda_{k_1})\sum_{\substack{k_2=1\\k_1\neq k_2\\|k_1-k_2|\leq a_T}}^{N/2}\Big|\frac{\sin(\frac{\pi (k_1-k_2)}{N}\Delta)}{\pi (k_1-k_2)}\Big|^2\nonumber\\
+&O(\frac{1}{a_T^{1-\tau}}+\frac{a_T}{N}+\frac{N^2}{a_T^3})\nonumber\\
=&\frac{2}{N}\sum_{k_1=a_T}^{N/2-a_T}\phi_{y_1,T}(j_1,\lambda_{k_1})\phi_{y_2,T}(j_2,\lambda_{k_1})f^2(\lambda_{k_1})\Big(\sum_{k_2=1}^{a_T}\frac{\sin^2(\pi k_2\frac{\Delta}{N})}{\pi^2k_2^2}\Big)+O(\frac{1}{a_T^{1-\tau}}+\frac{a_T}{N}+\frac{N^2}{a_T^3}) \nonumber \\
=& \frac{1}{N}\sum_{\substack{k_1=1}}^{N/2}\phi_{y_1,T}(j_1,\lambda_{k_1})\phi_{y_2,T}(j_2,\lambda_{k_1})f^2(\lambda_{k_1})\frac{(N-\Delta)\Delta}{N^2}+O(\frac{1}{a_T^{1-\tau}}+\frac{a_T}{N}+\frac{N^2}{a_T^3}) \nonumber,
\end{align}
where we used the identities $\sin^2(x)=\frac{1}{2}(1-\cos(2x))$ and $\sum_{k=1}^{\infty}\frac{\cos(kx)}{k^2}=\big(\frac{x-\pi}{2}\big)^2-\frac{\pi^2}{12}$ in the last step [see \cite{jolley}]. By combining the above results for $V_{k_1=k_2}$ and $V_{k_1\neq k_2}$ and proceeding completely analogously in the treatment of the second summand of $B_T(j_1,j_2)$ we obtain
\begin{align*}
B_T(j_1,j_2)=\big(1-\frac{\Delta}{N}\big)\frac{2}{N}\sum_{\substack{k=1}}^{N/2}\phi_{y_1,T}(j_1,\lambda_{k})\phi_{y_2,T}(j_2,\lambda_{k})f^2(\lambda_{k})
+O(\frac{1}{a_T^{1-\tau}}+\frac{a_T}{N}+\frac{N^2}{a_T^3}),
\end{align*}
uniformly with respect to $(j_1,j_2)$ if $\Delta=|j_1-j_2|<N$. For the case $\frac{1}{c}\leq v_i\leq1-\frac{1}{c}$ for $i=1,2$ the definition of the functions $\phi_{y_1,T}$ and $\phi_{y_2,T}$ now yields
\begin{align}
&N\Cov(\hat D_T(\phi_{y_1,T}),\hat D_T(\phi_{y_2,T}))=\sum_{\substack{(j_1,j_2)\in A_T(v_1,v_2)\\\Delta<N}}B_T(j_1,j_2)+O(\frac{1}{a_T^{1-\epsilon}}+\frac{a_T}{N}+\frac{N^2}{a_T^3})\nonumber\\
&\rightarrow \begin{cases}
0 \quad&\text{if }\frac{2}{c}\leq v_2-v_1\\
-[2-(v_2-v_1)c]\frac{1}{\pi}\int_0^{\min(\omega_1,\omega_2)\pi}f^2(\lambda)d\lambda \quad&\text{if } \frac{1}{c}<v_2-v_1\leq\frac{2}{c}\\
[2-3(v_2-v_1)c]\frac{1}{\pi}\int_0^{\min(\omega_1,\omega_2)\pi}f^2(\lambda)d\lambda \quad&\text{if } 0\leq v_2-v_1\leq\frac{1}{c},\\
\end{cases}\nonumber
\end{align}
and the assertion for the cases $v_i<\frac{1}{c}$ and $v_i>1-\frac{1}{c}$ for at least one $i\in\{1,2\}$ follows by similar arguments.
\medskip

For the proof of part (iii) we use the notation $Y_{i,1}:=Z_{j_i-   N/2   +p_i-m_i}$,  $Y_{i,2}:=Z_{j_i-   N/2   +q_i-n_i}$,
$$A_{T,l}(v_1,...,v_l):=\{\lfloor v_1T\rfloor-   N/2   ,\lfloor v_1T\rfloor+   N/2\}\times ... \times   \{\lfloor v_lT\rfloor-   N/2   ,\lfloor v_lT\rfloor+   N/2   \}$$
and obtain as in the proof of Theorem 2.1 (2) in \cite{detprevet2011b} that
\begin{align*}
\cum(\sqrt{N}\hat D_T(\phi_{y_1,T}),...,\sqrt{N}\hat D_T(\phi_{y_l,T}))=\sum_{\nu}V(\nu),
\end{align*}
where $V(\nu)$ is defined
by
\begin{align}
&V(\nu):=\frac{1}{N^{\frac{3l}{2}}(2\pi)^l}\sum_{(j_1,...,j_l) \in A_{T,l}(v_1,...,v_l)}\sum_{k_1,...,k_l=1}^{   N/2   }\sum_{m_1,...,m_l= 0}^{\infty}\sum_{n_1,...,n_l= 0}^{\infty}\sum_{p_1,...,p_l,q_1,...,q_l=0}^{N-1}\nonumber\\
&\phantom{...}\times\prod_{s=1}^l\phi_{y_s,T}(j_s,\lambda_{k_s})\psi_{m_s}\psi_{n_s}e^{-\lambda_{k_s}(p_s-q_s)}\text{cum}(Y_{a,b};(a,b)\in\nu_1)\cdots\text{cum}(Y_{a,b};(a,b)\in\nu_l),\nonumber
\end{align}
and the summation is carried out over all indecomposable partitions $\nu=(\nu_1,...,\nu_k)$ of the table
\begin{equation}\label{schema}
\begin{matrix}
Y_{1,1}&Y_{1,2}\\
\vdots&\vdots\\
Y_{l,1}&Y_{l,2}\\
\end{matrix}
\end{equation}
[see Theorem 2.3.2 in \cite{brillinger1981}]. Due to the Gaussianity of the $Y_{i,j}$ we only have to consider partitions $\nu=(\nu_1,...,\nu_l)$ with $l$ elements and without loss of generality we restrict ourselves  to the indecomposable partition
\begin{align}\label{nu}
\bar{\nu}:=\bigcup_{i=1}^{l-1}(Y_{i,1},Y_{i+1,2})\cup(Y_{l,1},Y_{1,2})
\end{align}
[all partitions which yield cumulants of the same order can be treated in the same way, while all other partitions yield cumulants of smaller order]. Simple calculations reveal
\begin{align}
V(\bar\nu)=&\frac{1}{N^{\frac{3}{2}l}(2\pi)^l}\sum_{(j_1,...,j_l) \in A_{T,l}(v_1,...,v_l)}\sum_{k_1,...,k_l=1}^{N/2}\prod_{s=1}^l\phi_{y_s,T}(j_s,\lambda_{k_s})\sum_{\substack{m_1,...,m_l,n_1,...,n_l=0}}^{\infty}\sum_{p_1,...,p_l}^{\star}\nonumber\\
\times& \prod_{s=1}^l\psi_{m_s}\psi_{n_s}e^{-i\lambda_{k_1}(p_1-p_l+m_l-n_1+j_1-j_{l})}\prod_{s=2}^le^{-i\lambda_{k_s}(p_s-p_{s-1}+m_{s-1}-n_s-j_{s-1}+j_s)},\label{cumm1}
\end{align}
where $\sum\limits^{\star}$  denotes the summation with respect to the conditions
\begin{align}
0\leq q_{i+1}&=p_i-m_i+n_{i+1}+j_i-j_{i+1}\leq N-1 \quad\text{ for }\quad i\in\{1,...,l-1\}\label{conditions1}\\
0\leq q_1&=p_l-m_l+n_1-j_1+j_{l}\leq N-1,\nonumber
\end{align}
[otherwise the cumulant vanishes due to the independence of the $Z_t$]. \eqref{conditions1} implies that
\begin{align}
&|n_{i+1}-m_i+j_i-j_{i+1}|\leq N \quad\text{ for }\quad i\in\{1,...,l-1\}\label{conditions2}\\
&|n_1-m_l-j_1+j_{l}|\leq N\nonumber
\end{align}
have to hold for each summand in \eqref{cumm1} not to vanish. If we combine this with the inequality
\begin{align*}
\frac{1}{N}\sum_{k=1}^{N/2}\phi_{y,T}(j,\lambda_k)e^{-i\lambda_kr}\leq\frac{C}{|r\text{ mod }N/2|},
\end{align*}
which holds uniformly with respect to $v$, $\omega$ and for all $r \in \N$ with $r$ mod $N/2\neq 0$ [see (A2) in \cite{eichler2008}], and the bound $\psi_l\leq Cl^{-1}$, which follows from \eqref{summ1}, we obtain that \eqref{cumm1} is bounded by (up to a constant)
\begin{align}
&\frac{1}{N^{l/2}}\sum_{(j_1,...,j_l) \in A_{T,l}(v_1,...,v_l)}\sum_{\substack{m_1...n_l=1\\|n_{i+1}-m_i+j_i-j_{i+1}\leq N\\|n_1-m_l-j_1+j_{l}\leq N}}^{\infty}\frac{1}{m_1\hdots m_l\cdot n_1\hdots n_l}\sum_{\substack{p_1,...,p_l=0\\|p_s-p_{s-1}+m_{s-1}-n_s-j_{s-1}+j_s|<N/2\\|p_1-p_l+m_l-n_1+j_1-j_l|<N/2}}^{N-1} \nonumber \\
\times&\prod_{s=1}^l\frac{1}{|p_s-p_{s-1}+m_{s-1}-n_s-j_{s-1}+j_s|}\prod_{s=1}^l1(p_s\notin \{z_{s1},z_{s2}\})\label{cumm2},
\end{align}
where we define $z_{s1}:=-p_{s-1}+m_{s-1}-n_s-j_{s-1}+j_s$ and $z_{s2}:=p_{s+1}+m_{s}-n_{s+1}-j_{s}+j_{s+1}$ (identifying $0$ with $l$ and $l+1$ with $1$) and use that the cases with $p_s=z_{s1}$, $p_s=z_{s2}$ or $|p_s-p_{s-1}+m_{s-1}-n_s-j_{s-1}+j_s|\geq N/2$ for some $s\in\{1,...,l\}$ are of the same or smaller order than \eqref{cumm2}. By proceeding completely analogously to the proof of Theorem 6.1 c)  in \cite{preussC} [with $d=0$] thereafter, it  follows that $V(\bar{\nu})$ is of order $o(1)$.  \\

\subsection{Proof of (\ref{cond0}) and (\ref{cond1}) in the proof of Theorem \ref{hauptsatz}}

Since the proof of (\ref{cond0}) is elementary we restrict ourselves to a treatment of \eqref{cond1}, and  by following the proof of Theorem 2.4 in \cite{dahlhaus1988}, it suffices to show
\be
|\cum_l \Big(\sqrt{N}(\hat{D}^{(1)}_T(y_1)-\hat{D}^{(1)}_T(y_2) \big))|\leq (2l)!C^ld_T(y_1,y_2)^l \label{alternativeaussage}
\ee
for all $(y_1,y_2)\in \mathcal{P}_T$, $l \in \N$. Without loss of generality we assume that $l$ is even and $\omega_1\leq\omega_2$. We use the notations $\phi(j,\lambda):=\phi^{(1)}_{y_1,T}(j,\lambda)-\phi^{(1)}_{y_2,T}(j,\lambda)$ and obtain as in the proof of part (1) (iii) that
\bea
\cum_l \Big(\sqrt{N} \big(\hat{D}^{(1)}_T(y_1)-\hat{D}^{(1)}_T(y_2) \big)\Big)=\sum_{\nu} V(\nu),
\eea
where $V(\nu)$ is defined as
\begin{align}
V(\nu&):=\frac{1}{N^{\frac{3l}{2}}(2\pi)^l}\sum_{j_1,...,j_l=1}^{T}\sum_{k_1,...,k_l=1}^{   N/2   }\sum_{m_1,...,m_l= 0}^{\infty}\sum_{n_1,...,n_l= 0}^{\infty}\sum_{p_1,...,p_l,q_1,...,q_l=0}^{N-1}\nonumber\\
&\times\prod_{s=1}^l\phi(j_s,\lambda_{k_s})\psi_{m_s}\psi_{n_s}e^{-\lambda_{k_s}(p_s-q_s)}\text{cum}(Y_{a,b};(a,b)\in\nu_1)\cdots\text{cum}(Y_{a,b};(a,b)\in\nu_l),\label{V1}
\end{align}
and the summation is carried out over all indecomposable partitions $\nu=(\nu_1,...,\nu_k)$ of \eqref{schema}. As in the proof of part (1) (iii),  we only have to consider partitions $\nu=(\nu_1,...,\nu_l)$ with $l$ elements due to the Gaussianity of the innovations and restrict ourselves  to the indecomposable partition of the form \eqref{nu}. The definition of the function $\phi(j,\lambda)$ implies that the summation over the set $\{1,...,T\}^l$ in \eqref{V1} needs only to be performed over the set
\begin{equation*}
A_T(v_1,v_2)^l:=\{\lfloor v_1T\rfloor+   N/2 ,\lfloor v_2T\rfloor+   N/2  \}^l.
\end{equation*}
By additionally exploiting that $\text{cum}(Z_i,Z_j)=\delta_{ij}$, we again obtain the $l$ relations \eqref{conditions1} between $p_i,q_i,m_i,n_i$ and  $j_i$, which have to be fulfilled such that the corresponding summands in \eqref{V1} do not vanish. We redefine the $m_i$ and $n_i$ [denoting $m_i$ by $n_i$ and replacing $n_i$ by $m_{i-1}$], which gives the following conditions
\begin{align}
0&\leq q_{i+1}=p_i+m_{i}-n_i+j_i-j_{i+1}\leq N-1\text{ for }i\in\{1,3,...,l-1\}\label{bed13.2}\\
0&\leq q_{i+1}=p_i+m_{i}-n_i+j_i-j_{i+1}\leq N-1\text{ for }i\in\{2,4,...,l-2\}\label{bed13.2.2}\\
0&\leq q_{1}=p_l+m_{l}-n_l+j_l-j_{1}\leq N-1\label{bed23.2}.
\end{align}
By employing the Cauchy Schwarz inequality we can bound $V(\bar{\nu})$ by $\sqrt{V_1(\bar{\nu})V_2(\bar{\nu})}$,
where
\begin{align}
V_1&(\bar{\nu}):=\frac{1}{N^{\frac{3l}{2}}(2\pi)^l}\sum_{(j_1,...,j_l)\in A_T(v_1,v_2)^l}\sum_{k_1,...,k_l=1}^{   N/2   }\Big(\prod_{s=1}^l\phi(j_s,\lambda_{k_s})\Big)\Big|\sum_{p_1,p_3...,p_{l-1}=0}^{N-1}\sum_{\substack{m_1,n_1,m_3,n_3,..\\.,m_{l-1},n_{l-1}= 0\\(\ref{bed13.2})}}^{\infty}\nonumber\\
\times&\prod_{s\in\{1,3,...,l-1\}}\psi_{m_s}\psi_{n_s}\prod_{s\in\{1,3,...,l-1\}}e^{-i(\lambda_{k_{s}}-\lambda_{k_{s+1}})p_s}\prod_{s\in\{1,3,...,l-1\}}e^{-i\lambda_{k_{s+1}}(n_s-m_{s}+j_{s+1}-j_s)}\Big|^2  \label{vquer2.1}
\end{align}
and $V_2(\bar{\nu})$ is defined analogously with the respective summations carried out over $p_i,m_i,n_i $ for even $i$ and with respect to conditions \eqref{bed13.2.2} and (\ref{bed23.2}). The term $V_1(\bar{\nu})$ is bounded by
\begin{align}
&\frac{C}{N^{3l/2}}\sum_{(j_1,...,j_l)\in A_T(v_1,v_2)^l}\sum_{k_1,k_2...,k_{l}=1}^{   N/2   }\Big(\prod_{s=1}^l\phi(j_s,\lambda_{k_s})\Big)\Big|H_T(j_1,j_2,...,j_l,\lambda_{k_1},\lambda_{k_2},...,\lambda_{k_{l}})\Big|\label{cum3},
\end{align}
where $H_T(j_1,j_2,...,j_l,\lambda_{k_1},\lambda_{k_2},...,\lambda_{k_{l}})$ equals
\bea
&&\sum_{\substack{p_1,p_3...,p_{l-1}\\\bar p_1,\bar p_3...,\bar p_{l-1}=0}}^{N-1}\sum_{\substack{m_1,n_1,m_3,n_3,...\\...,m_{l-1},n_{l-1}= 0\\(\ref{bed13.2})}}^{\infty}\sum_{\substack{\bar m_1,\bar n_1,\bar m_3,\bar n_3,...\\...,\bar m_{l-1},\bar n_{l-1}= 0\\ \overline{(\ref{bed13.2})}}}^{\infty}\prod_{s\in\{1,3,...,l-1\}}\Big( e^{-i\lambda_{k_s}(p_s-\bar p_s)}\psi_{m_s}\psi_{n_s}\psi_{\bar m_s}\psi_{\bar n_s}\\
&&\times e^{-i\lambda_{k_{s+1}}(-p_s+\bar p_s+n_s-m_s-\bar n_s+\bar m_s)}\Big)\nonumber
\eea
and $\overline{(\ref{bed13.2})}$ denotes the condition (\ref{bed13.2}) with $p_i,m_i,n_i$ replaced by $\bar p_i,\bar m_i,\bar n_i$. Simple calculations show that \eqref{cum3} is equal to
\begin{align}
&\frac{C}{N^{3l/2}}\sum_{(j_1,...,j_l)\in A_T(v_1,v_2)^l}\sum_{k_1,k_2...,k_{l}=1}^{   N/2   }\Big(\prod_{s=1}^l\phi(j_s,\lambda_{k_s})\Big)\sum_{\substack{m_1,n_1,m_3,n_3,...\\...,m_{l-1},n_{l-1}= 0\\(\ref{bed13.2})}}^{\infty}\sum_{\substack{\bar m_1,\bar n_1,\bar m_3,\bar n_3,...\\...,\bar m_{l-1},\bar n_{l-1}= 0\\ \overline{}}}^{\infty}\label{cum4}\\
\times&\prod_{s\in\{1,3,...,l-1\}}\psi_{m_s}\psi_{n_s}\psi_{\bar m_s}\psi_{\bar n_s}\prod_{s\in\{2,4,...,l\}}e^{-i\lambda_{k_s}(n_{s-1}-m_{s-1}  -\bar n_{s-1}+\bar m_{s-1})}\sum_{\substack{p_1,p_3...,p_{l-1}=0\\(\ref{bed13.2})}}^{N-1}\sum_{\substack{\bar p_1,\bar p_3...,\bar p_{l-1}=0\\ \overline{(\ref{bed13.2})}}}^{N-1}\nonumber\\
\times&\prod_{s\in\{1,3,...,l-1\}}e^{-i(\lambda_{k_s}-\lambda_{k_{s+1}})(p_s-\bar p_s)}\nonumber\\
=&\Big(\frac{C}{N^3}\sum_{(j_1,j_2)\in A_T(v_1,v_2)^2}\sum_{k_1,k_2=1}^{   N/2   }\phi(j_1,\lambda_{k_1})\phi(j_2,\lambda_{k_2})\sum_{m,n,\bar m,\bar n=0}^{\infty}\psi_m\psi_n\psi_{\bar m}\psi_{\bar n}e^{-i\lambda_{k_2}(n-m-\bar n+\bar m)}\nonumber\\
\times&\sum_{\substack{p_1,\bar p_1=0\\0\leq p_1+m_1-n_1+j_1-j_2\leq N_1\\0\leq \bar p_1+\bar m_1-\bar n_1+j_1-j_2\leq N_1}}^{N-1}e^{-i(\lambda_{k_1}-\lambda_{k_{2}})(p_1-\bar p_1)}\Big)^{l/2}=:\Big(\sum_{(j_1,j_2)\in A_T(v_1,v_2)^2}\tilde{B}_T(j_1,j_2)\Big)^{l/2}\nonumber,
\end{align}
where the quantities $\tilde{B}_T(j_1,j_2)$ are defined in an obvious manner. By the same arguments as provided in the calculation of the covariances and by the definition of $\phi$ one obtains that $\tilde{B}_T(j_1,j_2)=0$ if $|j_1-j_2|>N$ and that \eqref{cum4} is bounded by
\begin{align*}
C\Big[\int_{\omega_1\pi}^{\omega_2\pi}f^2(\lambda)d\lambda+\frac{2|\lfloor u(v_1,T)T\rfloor-\lfloor u(v_2,T)T\rfloor|}{N}\int_0^{\omega_1\pi}f^2(\lambda)d\lambda\Big]^{l/2}\leq Cd_T(y_1,y_2)^{l/2}.
\end{align*}
Since the same upper bound holds for $V_2(\bar \nu)$, the claim follows.

\subsection{Proof of Theorem \ref{hauptsatz2}}
As in the proof of Theorem \ref{hauptsatz}  we restrict ourselves to the case $d=1$ and emphazise that $C$ denotes a constant, which does not depend on the sample size and can vary from line to line in the calculations.\\

{\bf Proof of part a):} In order to show the claim it is sufficient to prove the following two properties [see. \cite{unifeconometrica}]:
\begin{itemize}
\item[(1)] For every $y:=(v,\omega)\in [0,1]$ we have
\begin{align}
\label{finitedimensional2}
N^\gamma \hat D_{T}(y) =o_P(1).
\end{align}
\item[(2)] For every $\eta, \epsilon>0$ there exists a $\delta>0$ such that
\begin{align}
\label{gleichstetigkeit2}
\lim_{T\rightarrow \infty} P\Big( \sup_{(y_1, y_2) \in [0,1]^2: d_2(y_1,y_2) <\delta } N^\gamma |\hat D_T(y_1)-\hat D_T(y_2)| > \eta \Big)< \epsilon ,
\end{align}
where $d_2(y_1,y_2)$ denotes the euclidean distance between $y_1=(v_1,\omega_1)$ and $y_2=(v_2,\omega_2)$.
\end{itemize}

The statement \eqref{finitedimensional2} follows by similar arguments as given in the proof of Theorem \ref{hauptsatz} a) and the details are omitted for the sake of brevity. For a proof of \eqref{gleichstetigkeit2} we employ the representation
\begin{align*}
\hat D_T(y)&=\frac{1}{N}\sum_{j=1}^{T}\sum_{k=1}^{   N/2   }\phi_{y,T}(j,\lambda_k)I_N(\frac{j}{T},\lambda_k)
\end{align*}
which was introduced in \eqref{REPRESENTATION}, where the functions $\phi_{y,T}$ for $y=(v,\omega)$ are defined in \eqref{PHI}. We define by
\be \label{rho2gamma}
\rho_{2,T,\gamma}(y_1,y_2):=\Big(\frac{1}{N^{2(1-\gamma)}}\sum_{j=1}^T\sum_{k=1}^N\big(\phi_{y_1,T}(j,\lambda_k)-\phi_{y_2,T}(j,\lambda_k)\big)^2\Big)^{1/2}
\ee
a semi-metric on the set
$$\mathcal{P}_T:=\{N/T,(N+1)/T,...,1-N/T\}\times\{1/N,...,1-1/N,1\},$$
and a straightforward calculation shows that  there exists a constant $A \in \R^+$ such that the inequality
\begin{align}
 \rho_{2,T,\gamma}(&y_1,y_2) \leq d_{T,\gamma}(y_1,y_2)\label{boundrho2}\\
:&=\frac{A/2}{N^{\frac{1}{2}-\gamma}}(\sqrt{|\omega_2-\omega_1|}+\sqrt{|v_2-v_1|})1_{\{\lfloor v_1T\rfloor=\lfloor v_2T\rfloor\}}+\frac{A}{N^{\frac{1}{2}-\gamma}}1_{\{\lfloor v_1T\rfloor\neq\lfloor v_2T\rfloor\}} \nonumber
\end{align}
holds. Note that $d_{T,\gamma}(\cdot,\cdot)$ defines a further semi-metric on $\mathcal{P}_T$, and for every $\gamma \in [0,1/2)$ we obtain for fixed $\delta>0$ and sufficiently large $T\in\N$
\begin{align}
\mathbb{P}(\sup\limits_{\substack{y_i\in[0,1]:\\ d_2(y_1,y_2)<\delta}}N^{\gamma}|\hat D_T(y_1)-\hat D_T(y_2) |>\eta)
&\leq &\mathbb{P}(\sup\limits_{\substack{y_i\in \mathcal{P}_T\\d_{T,\gamma}(y_1,y_2)<\delta}}N^{\gamma}|\hat D_T(y_1)-\hat D_T(y_2) |>\eta)\label{ungl32}.
\end{align}

So it suffices to prove that the probability on the right hand side of \eqref{ungl32} is arbitrary small for $T$ sufficiently large. This is a consequence of the following assertions:
\begin{itemize}
\item[i)]$\lim\limits_{\kappa \rightarrow 0}\lim\limits_{T \rightarrow \infty}J_{T,\gamma}(\kappa)=0$, where $J_{T,\gamma}$ denotes the covering integral of the set $\mathcal{P}_T$ with respect to the metric $d_{T,\gamma}$.
\item[ii)] There exists a constant $C \in \R^+$ such that for all $y_1,y_2\in [0,1]^2$ and integers $l\in\N$
\be
\label{cond12}
\Big|\cum_l \Big(N^\gamma (\hat{D}_T(y_1)-\hat{D}_T(y_2) \big)\Big)\Big|\leq (2l)!C^ld_{T,\gamma}(y_1,y_2)^k.
\ee
\end{itemize}

Proof of part i): The construction of $d_{T,\gamma}$ implies that there exists a constant $A>0$ such that  $C(u,d_{T,\gamma},\mathcal{P}_T)\leq TN$ for $u\leq A/N^{1/2-\gamma}$ and $C(u,d_{T,\gamma},\mathcal{P}_T)=1$ for $u\geq A/N^{1/2-\gamma}$. This yields that, for any $\kappa>0$, the quantity $J_{T,\gamma}(\kappa)$ is bounded by
\begin{align*}
&J_{T,\gamma}(A/N^{1/2-\gamma})+1_{\{\kappa>A/N^{1/2-\gamma}\}}\int_{A/N^{1/2-\gamma}}^{\kappa}\big[\log\big(\frac{48C(u,d_{T,\gamma},\mathcal{P}_T)^2}{u}\big)\big]^2du\\
\leq &\int_0^{A/N^{1/2-\gamma}}[\log(48T^2N^2)^2-2\log(48T^2N^2)\log(u)+\log(u)^2]du+\big|\int_{A/N^{1/2-\gamma}}^{\kappa}\big[\log\big(\frac{48}{u}\big)\big]^2du\big|,
\end{align*}
and it is straightforward to see that the right handside of the above expression converges to zero as $T \rightarrow \infty$ and $\kappa\rightarrow 0$.

Proof of part ii): We assume without loss of generality that $l$ is even and only consider the first summand $\{\hat D_T^{(1)}(y)\}_{y \in [0,1]^2}$. We then obtain as in the proof of \eqref{cond1} (by using an additional symmetry argument) that
\bea
\cum_l \Big(N^\gamma \big(\hat{D}_T^{(1)}(y_1)-\hat{D}_T^{(1)}(y_2) \big)\Big)=\sum_{\nu} V_\gamma(\nu),
\eea
where $V_\gamma(\nu)$ is defined as
\begin{align}
V_\gamma(\nu):=&\frac{1}{N^{l(2-\gamma)}(4\pi)^l}\sum_{j_1,...,j_l=1}^{T}\sum_{k_1,...,k_l=-\lfloor\frac{N-1}{2}\rfloor}^{   N/2   }\sum_{m_1,...,m_l= 0}^{\infty}\sum_{n_1,...,n_l= 0}^{\infty}\sum_{p_1,...,p_l,q_1,...,q_l=0}^{N-1}\prod_{s=1}^l\phi(j_s,\lambda_{k_s})\nonumber\\
&\times\prod_{s=1}^l[\psi_{m_s}\psi_{n_s}e^{-\lambda_{k_s}(p_s-q_s)}]\text{cum}(Y_{a,b};(a,b)\in\nu_1)\cdots\text{cum}(Y_{a,b};(a,b)\in\nu_l),\label{V12}
\end{align}
and the summation is carried out over all indecomposable partitions $\nu=(\nu_1,...,\nu_k)$ of \eqref{schema} and without loss of generality we restrict ourselves  to the indecomposable partition defined in \eqref{nu}.  By proceeding completely analogously to the proof of \eqref{cond1} we obtain
\begin{equation*}
V_\gamma(\bar{\nu})\leq\sqrt{V_{1,\gamma}(\bar{\nu})V_{2,\gamma}(\bar{\nu})},
\end{equation*}
where
\begin{align}
V_{1,\gamma}(\bar{\nu}):=&\frac{1}{N^{l(2-\gamma)}(4\pi)^l}\sum_{(j_1,...,j_l)\in A_T(v_1,v_2)^l}\sum_{k_1,...,k_l=-\lfloor\frac{N-1}{2}\rfloor}^{   N/2   }\prod_{s\in\{1,3,...,l-1\}}\phi^2(j_s,\lambda_{k_s})\nonumber\\
\times&\Big|\sum_{p_1,p_3...,p_{l-1}=0}^{N-1}\sum_{\substack{m_1,n_1,m_3,n_3,...,m_{l-1},n_{l-1}= 0\\(\ref{bed13.2})}}^{\infty}\prod_{s\in\{1,3,...,l-1\}}\psi_{m_s}\psi_{n_s}\prod_{s\in\{1,3,...,l-1\}}e^{-i(\lambda_{k_{s}}-\lambda_{k_{s+1}})p_s}\nonumber\\
\times&\prod_{s\in\{1,3,...,l-1\}}e^{-i\lambda_{k_{s+1}}(n_s-m_{s}+j_{s+1}-j_s)}\Big|^2\label{vquer2.2}
\end{align}
and $V_{2,\gamma}(\bar{\nu})$ is defined analogously with the respective summations carried out over $j_i,p_i,m_i,n_i $ for even $i$ and with respect to condition (\ref{bed23.2}). The term $V_{1,\gamma}(\bar{\nu})$ is bounded by
\begin{align}
&\frac{C}{N^{l(1-\gamma)}}\sum_{\substack{(j_1,...,j_l)\\ \in A_T(v_1,v_2)^l}}\sum_{k_1,k_3...,k_{l-1}=-\lfloor\frac{N-1}{2}\rfloor}^{   N/2   }\prod_{s\in\{1,3,...,l-1\}}\phi^2(j_s,\lambda_{k_s}) |H_{T,2}(j_1,j_2,...,j_l,\lambda_{k_1},\lambda_{k_3},...,\lambda_{k_{l-1}})|\nonumber,
\end{align}
where $H_{T,2}(j_1,j_2,...,j_l,\lambda_{k_1},\lambda_{k_3},...,\lambda_{k_{l-1}})$ equals
\bea
&&\frac{1}{N^{l}}\sum_{\substack{p_1,p_3...,p_{l-1}\\\bar p_1,\bar p_3...,\bar p_{l-1}=0}}^{N-1}\sum_{\substack{m_1,n_1,m_3,n_3,...\\...,m_{l-1},n_{l-1}= 0\\(\ref{bed13.2})}}^{\infty}\sum_{\substack{\bar m_1,\bar n_1,\bar m_3,\bar n_3,...\\...,\bar m_{l-1},\bar n_{l-1}= 0\\ \overline{(\ref{bed13.2})}}}^{\infty}\prod_{s\in\{1,3,...,l-1\}}e^{-i\lambda_{k_s}(p_s-\bar p_s)}\prod_{s\in\{1,3,...,l-1\}}\psi_{m_s}\psi_{n_s}\psi_{\bar m_s}\psi_{\bar n_s}\nonumber\\
&&\times\sum_{k_2,...,k_l=-\lfloor\frac{N-1}{2}\rfloor}^{   N/2   }\prod_{s\in\{1,3,...,l-1\}}e^{-i\lambda_{k_{s+1}}(-p_s+\bar p_s+n_s-m_s-\bar n_s+\bar m_s)}
\eea
and $\overline{(\ref{bed13.2})}$ denotes the condition (\ref{bed13.2}) with $p_i,m_i,n_i$ replaced by $\bar p_i,\bar m_i,\bar n_i$. The well know identity
\begin{align*}
\frac{1}{N}\sum_{k=-\lfloor\frac{N-1}{2}\rfloor}^{   N/2   }\exp(\lambda_kp)=
\begin{cases}
1&\text{ if }p=l N \text{ for some } l\in\mathbb{N}\\
0 &\text{ else}
\end{cases}
\end{align*}
yields that for a fixed choice of $p_i,m_i,n_i,\bar m_i,\bar n_i$ there is only one choice for $\bar p_i$ such that the corresponding summand does not vanish. Thus we obtain
\begin{align*}
&\sum_{(j_2,j_4...,j_{l})\in A_T(v_1,v_2)^{l/2}} |H_{T,2}(j_1,j_2,...,j_l,\lambda_{k_1},\lambda_{k_3},...,\lambda_{k_{l-1}})|\\
\leq&\sum_{(j_2,j_4...,j_{l})\in A_T(v_1,v_2)^{l/2}}\sum_{\substack{m_1,m_3...,m_{l-1}\\ \bar m_1,\bar m_3...,\bar m_{l-1}\\n_1,n_3,...,n_{l-1}\\ \bar n_1,\bar n_3,...,\bar n_{l-1}= 0}}^{\infty}\prod_{s\in\{1,3,...,l-1\}}|\psi_{m_s}\psi_{n_s}\psi_{\bar m_s}\psi_{\bar n_s}|=4^{l/2}\Big(\sum_{n= 0}^{\infty}|\psi_n|\Big)^{2l},
\end{align*}

which is bounded by $C^l$ for some $C>0$ due to \eqref{summ1}. Observing \eqref{vquer2.2} this implies that $V_{1,\gamma}(\bar{\nu})$ is bounded by
\bea
&&\frac{C^l}{N^{l(1-\gamma)}}\sum_{j_1,j_3...,j_{l-1}=1}^{T}\sum_{k_1,k_3...,k_{l-1}=-\lfloor\frac{N-1}{2}\rfloor}^{   N/2   }\Big(\phi^2(j_1,\lambda_{k_1})\phi^2(j_3,\lambda_{k_3})...\phi^2(j_{l-1},\lambda_{k_{l-1}})\Big)
\eea

which equals $C^l\rho_{2,T,\gamma}(y_1,y_2)^l$. Since the same upper bound can be obtained for $V_2(\bar \nu)$, \eqref{cond12} follows with \eqref{boundrho2}. \\

{\bf Proof of part b):} Assertion b) of Theorem \ref{hauptsatz2} can be established using the same arguments as in the proof of part b) of Theorem \ref{hauptsatz}.$\hfill \Box$

\subsection{Proof of Theorem \ref{theorembootstrap2}}
We again consider only the case $d=1$ and suppress the argument $T$, when referring to the sequence $p=p(T)$. Note that due to Assumption \ref{annahmenbootstrap} the stationary process $X_t$ with spectral density $g(\lambda)=\int_0^1 f(u,\lambda) du$ possesses an MA($\infty$) representation
\be
\label{XtMAunendlich}
X_t=\sum_{l=0}^\infty \psi_l Z_{t-l}.
\ee
Moreover, under the null hypothesis we have $g\equiv f$ and the functions $\psi_l(u)$ in \eqref{defX_t} satisfy for all $u\in[0,1]$ $\psi_l(u)=\psi_l$ ($l\in\N$). Let $\{X_t^*\}_{t\in\Z}$ denote the process which is obtained by replacing $Z_t$ in \eqref{XtMAunendlich} by the bootstrap replicates $Z_t^*$ from step 2) in Algorithm \ref{algo1}, then $\{X_t\}_{t\in\Z}$ and $\{X_t^*\}_{t\in\Z}$ obviously have the same distribution. Thus, if we define $\hat D_{T,a}^*(y)$ as $\hat D_T(y)$ where the random variables $X_{t,T}$ are replaced by $X_t^*$, then part a) of Theorem \ref{theorembootstrap2} follows directly.

From the arguments presented in the proof of Theorem \ref{hauptsatz} we obtain that there exist $\tilde y\in[0,1]^2$ such that $|\hat D_{T,a}^*(\tilde y)|^2 \geq L/N$ for some constant $L \in \R^+$. So it suffices to prove
\begin{equation}\label{zz2}
\sqrt{N}\sup\limits_{(y)\in[0,1]^2}|\hat D_T^*(y)-\hat D_{T,a}^*(y)|=o_P(1),
\end{equation}
which is a consequence of the following two statements [see the proof of Theorem \ref{hauptsatz2}]:
\begin{itemize}
\item[(1)] For every $y=(v,\omega)\in [0,1]$ we have $\sqrt{N}(\hat D_T^*(y)-\hat D_{T,a}^*(y)) =o_P(1)$.
\item[(2)] For every $\eta, \epsilon>0$ there exists a $\delta>0$ such that
\begin{align}
\label{gleichstetigkeitbootstrap}
\lim_{T\rightarrow \infty} P\Big( \sup_{\substack{y_1, y_2 \in [0,1]^2:\\ d_2(y_1,y_2) <\delta }} \sqrt{N} |(\hat D_T^*(y_1)-\hat D_{T,a}^*(y_1))-(\hat D_T^*(y_2)-\hat D_{T,a}^*(y_2))| > \eta \Big)< \epsilon.
\end{align}
\end{itemize}

For the sake of brevity we restrict ourselves to a proof of \eqref{gleichstetigkeitbootstrap} since it is by far the more complicated one.
For this purpose we introduce the process
\be
X_{t}^{AR}(p):=\sum_{j=1}^pa_{j,p}X_{t-j}^{AR}(p)+Z_t^{AR}(p), \label{arptrue}
\ee
where $a_{1,p},...,a_{p,p}$ are defined in \eqref{arpfitt} and the $Z_j^{AR}(p)$ are independent centered Gaussian random variables with variance $\sigma_p=\E\big[\big(X_t-\sum_{j=1}^pa_{j,p}X_{t-j}\big)^2\big]$. Because $\{X_t\}_{t\in\mathbb{Z}}$ is the stationary process with spectral density $g(\lambda)=\int_0^1f(u,\lambda) du$, the process $\{X_{t}^{AR}(p)\}_{t\in\mathbb{Z}}$ corresponds to the 'best' approximation of the process \eqref{Ar1} by an AR($p$) model. It now follows from Lemma 2.3 in \cite{kreisspappol2011} that for sufficiently large $T$ the approximating process $\{X_t^{AR}(p)\}_{t\in\Z}$ and the bootstrap analogon $X_{t,T}^*$ in \eqref{fittedarp} have MA($\infty$) representations, i.e.
\begin{equation}\label{statX2}
X_t^{AR}(p)=\sum_{j=0}^{\infty}\psi_{j}^{AR}(p)Z^{AR}_{t-j}(p) \quad \text{and} \quad X_{t,T}^*=\sum_{j=0}^{\infty}\hat{\psi}_{j}^{AR}(p)Z^*_{t-j}.
\end{equation}
We define $W_T$ as the event
\be
\label{upperboundcond}
\sqrt{N}\sum_{m,n=0}^{\infty}|\psi_{m}\psi_{n}-\hat\psi^{AR}_{m}(p)\hat\psi^{AR}_{n}(p)| \leq 1
\ee
and $\bar D_T(y):=\sqrt{N}(\hat D_T^*(y)-\hat D_{T,a}^*(y)) \times 1_{W_T}$. The claim \eqref{gleichstetigkeitbootstrap} then follows  as  in the proof of \eqref{gleichstetigkeit2} from the following two assertions:
\begin{itemize}
\item[i)] $\mathbb{P}(W_T)\rightarrow 1$.
\item[ii)] There exists a constant $C \in \R^+$ such that for all $y_1=(v_1,\omega_1)$ and $y_2=(v_2,\omega_2)$
\begin{align}\label{ineqcum}
|\text{cum}_l(\bar D_T(y_1)-\bar D_T(y_2))|\leq (2l)!C^l\rho_{2,T,0}(y_1,y_2)^l, \quad \forall l \in \N,
\end{align}
where $\rho_{2,T,0}$ was defined in \eqref{rho2gamma}.
\end{itemize}

We will at first comment on the second statement \eqref{ineqcum} and justify the first part afterwards. If we define $V_{1,\gamma}(\bar \nu)^*$ as $V_{1,\gamma}(\bar \nu)$ in \eqref{vquer2.1}, with the quantities
\be
\label{maohnebootstrap}
\prod_{s\in\{1,3,...,l-1\}}\psi_{m_s}\psi_{n_s}
\ee
being replaced by
\be
\label{mamitbootstrap}
\E \Big(\prod_{s\in\{1,3,...,l-1\}}\big(\psi_{m_s}\psi_{n_s}-\hat \psi_{m_s}^{AR} \hat\psi_{n_s}^{AR}\big)\Big),
\ee
then we obtain by the same arguments as given in the proof of Theorem \ref{hauptsatz2} that $V_{1,\gamma}(\bar \nu)^*$ is bounded by
\bea
\frac{1}{(4\pi)^lN^{l/2}}\sum_{\substack{(j_1,...,j_l)\\\in A_T(v_1,v_2)}}\sum_{k_1,k_3...,k_{l-1}=-\lfloor\frac{N-1}{2}\rfloor}^{   N/2   }\prod_{s\in\{1,3,...,l-1\}}\phi^2(j_s,\lambda_{k_s}) |H_{T,2}^*(j_1,j_2,...,j_l,\lambda_{k_1},\lambda_{k_3},...,\lambda_{k_{l-1}})|,\label{z3}
\eea

where $H_{T,2}^*(j_1,j_2,...,j_l,\lambda_{k_1},\lambda_{k_3},...,\lambda_{k_{l-1}})$ is the corresponding bootstrap analogue of \linebreak $H_{T,2}(j_1,j_2,...,j_l,\lambda_{k_1},\lambda_{k_3},...,\lambda_{k_{l-1}})$, i.e. the quantity  \eqref{maohnebootstrap} is replaced by \eqref{mamitbootstrap}. Combining the arguments in the treatment of $H_{T,2}(j_1,...,\lambda_{k_{l-1}})$ with \eqref{upperboundcond} we obtain
\bea
\sum_{(j_2,j_4,...,j_l)\in A_T(v_1,v_2)^{l/2}}N^{l/2}|H_{T,2}^*(j_1,j_2,...,j_l,\lambda_{k_1},\lambda_{k_3},...,\lambda_{k_{l-1}})|\leq4^{l/2}
\eea
on the set $W_T$. This yields \eqref{ineqcum} as in the proof of \eqref{gleichstetigkeit2}, and it remains to prove $\mathbb{P}(W_T)\rightarrow 1$. Because of
\begin{align*}
\sum_{m,n=0}^{\infty}|\psi_{m}\psi_{n}-\hat\psi^{AR}_{m}(p)\hat\psi^{AR}_{n}(p)|
\leq \sum_{m=0}^{\infty}|\psi_m|\sum_{n=0}^{\infty}|\psi_n-\hat\psi^{AR}_{n}(p)|+\sum_{m=0}^{\infty}|\hat\psi^{AR}_m(p)|\sum_{n=0}^{\infty}|\psi_n-\hat\psi^{AR}_{n}(p)|
\end{align*}
it is sufficient to show
\begin{align}\label{zuzeigen}
\sqrt{N}\sum_{n=0}^{\infty}|\psi_n-\hat\psi^{AR}_{n}(p)|=o_P(1)\quad\text{ and }\quad\sum_{m=0}^{\infty}|\hat\psi^{AR}_m(p)|=O_P(1).
\end{align}

We now use Lemma 2.3 of \cite{kreisspappol2011} which implies that the polynomials $A_{p}(z):=1-\sum_{k=1}^{p}a_{j,p}z^k$ and $\hat A_{p}(z):=1-\sum_{k=1}^{p}\hat a_{j,p}z^k$ have no roots within the closed unit disc $\{z|\text{ }|z|\leq1+\frac{1}{p}\}$ if $p$ is sufficiently large.  An application of Cauchy´┐¢s inequality for holomorphic functions to the difference
$$\hat A^{-1}_{p}(z)-A^{-1}_{p}(z):=\sum_{k=1}^{\infty}[\hat\psi_{k}^{AR}(p)-\psi_{k}^{AR}(p)]z^k,$$
as in the proof of Lemma 2.5 in \cite{kreisspappol2011}, yields
\begin{align*}
|\hat{\psi}_l^{AR}(p)-\psi_l^{AR}(p)|&\leq\Big(\frac{1}{1+1/p}\Big)^l\max\limits_{|z|=1+\frac{1}{p}}\frac{\big|A_{p}(z)-\hat A_{p}(z)\big|}{\big|A_{p}(z)\hat A_{p}(z)\big|}\\
&\leq\Big(\frac{1}{1+1/p}\Big)^l\max\limits_{|z|=1+\frac{1}{p}}\frac{\sum_{k=1}^{p}|\hat a_{k,p}-a_{k,p}|(1+\frac{1}{p})^k}{\big|A_{p}(z)\hat A_{p}(z)\big|}\\
&=p\Big(1+\frac{1}{p}\Big)^{-l}O_P(\sqrt{\log(T)/T})
\end{align*}
uniformly with respect to $l$ and $p$. This bound implies
\begin{equation} \label{teil1}
\sqrt{N}\sum_{n=0}^{\infty}|\hat\psi^{AR}_{n}(p)-\psi^{AR}_{n}(p)|=O_P\left(\frac{p^2_{\text{max}}(T)\sqrt{\log(T)} \sqrt{N}}{\sqrt{T}}\right)=o_P(1),
\end{equation}

and Lemma 2.4 of \cite{kreisspappol2011}\ yields that (for sufficiently large $T$)
\begin{align}\label{teil2}
\sqrt{N}\sum_{n=0}^{\infty}|\psi^{AR}_{n}(p)-\psi_n|\leq\frac{\sqrt{N}}{p}\sum_{j=p+1}^{\infty}j|a_j|=O_P\left(\frac{\sqrt{N}}{p_{min}(T)}\right)=o_P(1).
\end{align}

\eqref{teil1} and \eqref{teil2} now imply the first part of (\ref{zuzeigen}). The second part can be easily derived from the first part and \eqref{summ1}. This yields part b). Part c) is of course the bootstrap analogue of Theorem \ref{hauptsatz2} a), and is therefore shown by combining the arguments given here with the reasoning in the proof of Theorem \ref{hauptsatz2} a).$\hfill \Box$

\end{document}